\documentclass[10pt]{amsart}

\newcommand{\smn}{S_{m,n}}
\newcommand{\Smn}{S_{m+M,n+N}}
\newcommand{\Dom}{\text{Dom}_{m+M,n+N}}
\newcommand{\dom}{\text{Dom}_{m,n}}
\newcommand{\ls}{\left<}
\newcommand{\rs}{\right>}
\newcommand{\lb}{\left[\left[}
\newcommand{\rb}{\right]\right]_s}
\newcommand{\pb}{\bar{p}}
\newcommand{\qb}{\bar{q}}
\newcommand{\mfp}{\mathfrak{p}}
\newcommand{\mfm}{\mathfrak{m}}

\newcommand{\mfq}{\mathfrak{q}}

\newcommand{\kdom}{K'\text{-Dom}_{m,n}}
\newcommand{\bkdom}{\bar{K}\text{-Dom}_{m,n}}
\newcommand{\ideal}{\mathcal{I}}
\newcommand{\xr}{[x/\rho]}
\newcommand{\bag}{G}
\newcommand{\baf}{F}

\theoremstyle{plain}
\newtheorem{lem}{Lemma}[section]
\newtheorem{thm}[lem]{Theorem}
\newtheorem{cor}[lem]{Corollary}

\theoremstyle{definition}
\newtheorem{defn}[lem]{Definition}

\newtheorem{rem}[lem]{Remark}
\newtheorem{ex}[lem]{Example}

\begin{document}

\author{Y. F\i rat \c{C}el\.{i}kler}
\title[Parameterized stratification and piece number]{Parameterized stratification and piece number of $D$-semianalytic sets}
\address{Mathematics Department\\ Purdue University\\ West Lafayette\\ IN 47907, USA}
\email{ycelikle@math.purdue.edu}
\subjclass[2000]{Primary (32B20); Secondary (03C10), (12J25)}

\begin{abstract}
We obtain results on the geometry of $D$-semianalytic and subanalytic subsets over a complete, non-trivially valued non-Archimedean field $K$, which is not necessarily algebraically closed. Among the results are a parameterized smooth stratification theorem and several results concerning the dimension of the $D$-semianalytic and subanalytic sets. We also extend Bartenwerfer's definition of piece number for analytic $K$-varieties to $D$-semianalytic sets and prove the existence of a uniform bound for the piece number of the fibers of a $D$-semianalytic set. We also establish a connection between the piece number and complexity of $D$-semianalytic sets which are subsets of the line and thereby give a simpler proof of the Complexity Theorem of Lipshitz and Robinson. We finish by proving that for each $D$-semianalytic $X$, there is a semialgebraic $Y$ such that one dimensional fibers of $X$ are among the one dimensional fibers of $Y$. This is an analogue of a theorem by van den Dries, Haskell and Macpherson.
\end{abstract}

\thanks{The author wishes to thank his advisor Professor Leonard Lipshitz for his support and helpful suggestions.}

\maketitle

\section{Introduction}

In this paper we continue the investigation which started in \cite{dimtheory} of basic geometric properties of $D$-semianalytic subsets of $(K^\circ)^m$, where $K$ is an {\em arbitrary} ({\em i.e.} not necessarily algebraically closed) non-trivially valued, complete non-Archimedean field and $K^\circ$ is its valuation ring. For most of our results we impose the additional condition {\em Char $K=0$}. 

One of our main results states the existence of a uniform bound on the piece numbers of fibers of $D$-semianalytic sets (Theorem \ref{boundedpn}). Next we establish a connection between the piece number and the complexity of one dimensional $D$-semianalytic sets to give a simpler proof of the existence of a uniform bound for the complexity of the one-dimensional fibers of $D$-semianalytic sets (Theorem \ref{complexity}) which was first proved by Lipshitz and Robinson in \cite{complex}. Along the way we prove a parameterized version of the Smooth Stratification Theorem for $D$-semianalytic sets (Theorem \ref{parass}) and new results on the dimension theory of $D$-semianalytic and subanalytic sets. Namely we show that the dimension of the boundary of a $D$-semianalytic set is less than the dimension of the $D$-semianalytic set itself (Theorem \ref{bndry}) and that the dimension of a finite union of subanalytic sets is equal to the supremum of the dimensions of the individual sets (Theorem \ref{densesub}). We finish by showing that for each $D$-semianalytic set $X$ there is a semialgebraic set $Y$ such that the one dimensional fibers of $X$ are among the one-dimensional fibers of $Y$ (Theorem \ref{thmA}). This is an analogue of Theorem A of \cite{onedim} by van den Dries, Haskell and Macpherson.

$D$-semianalytic subsets of $(K^\circ)^{m}$ are finite unions of sets of points which satisfy finitely many norm inequalities (both strict and non-strict) between $D$-functions. These functions are obtained from members of rings of separated power series $\smn$ through the use of the restricted division operators $D_0$ and $D_1$, and substitution (see Definitions \ref{Dop} and \ref{Dfunc}). The rings of separated power series $\smn$ are superrings of the Tate algebras $T_{m+n}$ and were first introduced by Lipshitz in \cite{isop} in order to handle the case where $K$ is non-discretely valued while proving that when {\em Char $K=0$} or {\em Char $K=p>0$} and $[K:K^p]<\infty$, given a coordinate projection $\pi$, for each affinoid analytic $K$-variety $X$, there is a bound $\Gamma$ such that the fibers of $X$ under $\pi$ have at most $\Gamma$ isolated points. In \cite{bart} Bartenwerfer introduced the notion of the {\em piece number} of an analytic $K$-variety and proved that under the same conditions on the characteristic of $K$ as in \cite{isop}, given a projection $\pi$ and an analytic $K$-variety $X$ which is the zeroset of an ideal of $\smn$, there is a bound $\Gamma$ such that the fibers of $\pi$ have piece numbers less than $\Gamma$. As the piece number dominates the number of isolated points, this was a strengthening of the main theorem of \cite{isop}.

In \cite{rigid} Lipshitz showed that the class of subanalytic sets ({\em i.e.} the projections of $D$-semianalytic sets) coincides with the class of the $D$-semianalytic sets when $K$ is algebraically closed thus obtaining a quantifier elimination theorem for the analytic theory of such fields and later this result was generalized to a more general class of {\em definable sets} by Lipshitz and Robinson in \cite{modcomp}. Meanwhile, in \cite{strat}, they obtained results on the geometry of subanalytic subsets of $(\bar{K}^\circ)^m$ defined by $D$-functions with coefficients from $K$ where $K$ is an arbitrary complete non-Archimedean field and $\bar{K}$ is an algebraically closed complete extension of $K$. Among those results were the Smooth Stratification Theorem and results on the dimension theory of subanalytic sets. It was a natural question to ask how much of these results could be carried over to $D$-semianalytic subsets of $(K^\circ)^m$ where $K$ is an {\em arbitrary} non-Archimedean complete field. Results similar to those in \cite{strat} and further results on the dimension theory of $D$-semianalytic sets in this more general context were obtained in \cite{dimtheory}. One of the main tools in obtaining these results, the Parameterized Normalization Lemma for $D$-semianalytic Sets of \cite{dimtheory}, plays an essential role in proving many of the results in this paper also. The proof of this theorem relies on the ideas and results of \cite{strat}, where much of the commutative algebra of the rings $\smn$ and quasi-affinoid algebras are worked out.

On the other hand a special type of $D$-semianalytic sets over an algebraically closed field $K$, the $R$-domains, stands out in the study of the quasi-affinoid geometry. These sets generalize the {\em rational domains} of affinoid geometry as defined in Definition 5 of subsection 7.2.3 of \cite{bgr} and some of their geometric properties were established by Lipshitz and Robinson in \cite{linepln} and \cite{complex}. One of the main results of \cite{linepln} is that such subsets of $K^\circ$ can be written as boolean combinations of discs (see Definition \ref{disc}). In \cite{complex} the complexity of an $R$-subdomain of $K^\circ$ was defined to be the number of discs appearing in such a combination and it was shown that given a $D$-semianalytic subset $X$ of $(K^\circ)^{m+1}$ there is a bound $\Gamma$ such that for each parameter $\pb\in (K^\circ)^{m}$, the fiber $X(\pb)$ of $\pb$ in $X$ differs from an $R$-domain of complexity less than $\Gamma$ by at most $\Gamma$ points. This result is analogous to the Theorem A of \cite{onedim} by van den Dries, Haskell and Macpherson which states that the one-dimensional fibers of a subanalytic (in the language of \cite{denefvdd}) subset of $\mathbb{Z}_p^{m+1}$ can be obtained as one-dimensional fibers of a semialgebraic (in Macintyre's Language) subset of some $\mathbb{Z}_p^{M+1}$. As explained in \cite{onedim}, the bound on ``complexity'' of one dimensional fibers follows immediately from this theorem.

The outline of this paper is follows. After preliminary definitions in Section 2, we will revisit the dimension theory of $D$-semianalytic sets in Section 3, proving new results on the dimension of the boundary of a $D$-semianalytic set. We will also prove several properties of the parameterized normalization process to establish the facts needed in the later sections. The main theorem of Section 4 is the Parameterized Smooth Stratification Theorem which lays out the groundwork for the the results of Section 5. Also in this section, a theorem on the dimension theory of subanalytic sets is proved through a similar argument. We start Section 5 by extending the definition of piece number for analytic $K$-varieties due to Bartenwerfer, to $D$-semianalytic sets and prove the piece number theorem for $D$-semianalytic sets. In Section 6, we will turn our attention to one-dimensional fibers of $D$-semianalytic sets and establish a relation between the piece number and the complexity of such fibers resulting in a new, simpler proof of the Complexity Theorem of \cite{complex}. This theorem is analogous to the Theorem A of \cite{onedim} in the $p$-adic setting, yet we will finish by providing another, more readily recognizable, analogue of Theorem A in our setting. We would like to note that one could obtain the results stating the existence of uniform bounds like the Complexity Theorem by using non-standard models approach of \cite{onedim}.

\section{Preliminaries}

In this section we will give the definitions and basic properties of our objects of study. In doing so we will follow \cite{RoSPS} by Lipshitz and Robinson where they were studied extensively. These objects are mainly the rings of separated power series and geometric and algebraic objects related to them. 

 $K$ denotes an arbitrary non-trivially valued non-Archimedean complete field.

\begin{defn}\label{separ}
Let $x=(x_1,...,x_m)$ and $\rho=(\rho_1,...,\rho_n)$ denote variables, fix a complete, quasi-Noetherian subring $E$ of $K^{\circ}$ (which also has to be a discrete valuation ring in case $Char\: K=p>0$) and let $\{a_i\}_{i\in\mathbb{N}}$ be a zero sequence in $K^\circ$, and $B$ be the  local quasi-Noetherian ring
$$(E[a_0,a_1,...]_{\{a\in E[a_0,a_1,...]:|a|=1\}})^{\wedge},$$
where $\text{}^\wedge$ denotes the completion in $|\cdot|$. Let $\mathcal{B}$ be the family of all such rings. Define the {\em separated power series ring} over $(E,K)$ of $(m,n)$ variables to be
$$\smn(E,K):=K\otimes_{K^\circ}\left(\lim_{\overrightarrow{B\in\mathcal{B}}}B\ls x\rs[[\rho]]\right).$$

For $f=\sum_{\alpha,\beta}a_{\alpha,\beta}x^{\alpha}\rho^{\beta}\in\smn(E,K)$, the {\em Gauss norm} of $f$ is defined as 
$$||f||:=\sup_{\alpha,\beta}|a_{\alpha,\beta}|.$$
\end{defn}

We will  write $\smn$ instead of $\smn(E,K)$ when the ring $E$ and the field $K$ are clear from the context. We will make use of two key results on algebra of these rings from \cite{RoSPS} throughout. Namely, these rings are Noetherian and we have suitable Weierstrass Preparation and Division Theorems over these rings.

The ring $\smn$ contains our ``global'' analytic functions. In other words, members of $\smn(E,K)$ are convergent and have a global power series expansion over the set $(K^\circ)^m\times(K^{\circ\circ})^n$ where $K^{\circ\circ}$ denotes the maximal ideal of $K^\circ$. Nevertheless we are also interested in more general class of analytic functions which may have different power series expansions at different localities. More precisely we wish to be able to take quotients of members of rings of separated power series as well as substitute them in other members. For this to work we need two restricted division operators as defined by Lipshitz in \cite{rigid}.

\begin{defn}\label{Dop}
The restricted division operator $D_0:(K^\circ)^2\rightarrow K^\circ$ is defined as 
$$D_0(x,y):=\begin{cases}x/y &\text{ if } |x|\leq|y|\neq 0\\0&\text{ otherwise.}\end{cases}$$
Whereas the restricted division operator $D_1:(K^\circ)^2\rightarrow K^{\circ\circ}$ is defined as
$$D_1(x,y):=\begin{cases}x/y &\text{ if } |x|<|y|\\0&\text{ otherwise.}\end{cases}$$
\end{defn}

Now we can define the $D$-functions which were used in defining $D$-semianalytic sets. 

\begin{defn}\label{Dfunc}
$D$-functions over $(K^\circ)^m\times(K^{\circ\circ})^n$ are inductively defined as follows. 

i) Any member $f$ of $\smn$ is a $D$-function over $(K^\circ)^m\times(K^{\circ\circ})^n$. 

ii) If $f,g$ are analytic $D$ functions over $(K^\circ)^m\times(K^{\circ\circ})^n$ and $h\in S_{m+1,n}$ (or  $S_{m,n+1}$) then $h(x,D_0(f,g),\rho)$ (or $h(x,\rho,D_1(f,g))$) is a $D$-function over $(K^\circ)^m\times(K^{\circ\circ})^n$.
\end{defn}

Note that what we call a $D$-function coincides with what was called an $\mathcal{L}_{\text{an}}^D$-term in \cite{rigid} and in \cite{modcomp}. However we will rarely use these terms in the rest of this paper and instead work with ``generalized rings of fractions over $\smn$''. Before we define such rings we need some more notation. 

A ring $B$ which is of the form $\smn/I$ for some ring of separated power series $\smn$ and ideal $I\subset\smn$ is called a {\em quasi-affinoid algebra}. Let $f,g\in B$,  $z$ be a variable not appearing in $\smn$ and $f,g\in\smn$ be two elements whose canonical images in $B$ are $\bar{f}$ and $\bar{g}$ respectively, then we will write $B\ls \bar{f}/\bar{g}\rs$ (or $B\ls z\rs/(\bar{g}z-\bar{f})$) and $B\lb \bar{f}/\bar{g}\rb$ (or $B\lb z\rb/(\bar{g}z-\bar{f})$) for the rings $S_{m+1,n}/(I\cdot S_{m+1,n}\cup \{gx_{m+1}-f\})$ and $S_{m,n+1}/(I\cdot S_{m,n+1}\cup\{g\rho_{n+1}-f\})$ respectively. More generally if $A$ is a quasi affinoid algebra, $y$ and $\lambda$ are multi-variables not appearing in the presentation of $A$ and 
$$\begin{array}{rcl}
B_1&=&A\ls y_1,...,y_{M_1}\rs\lb\lambda_1,...,\lambda_{N_1}\rb/I_1\\
B_2&=&A\ls y_{M_1+1},...,y_{M_2}\rs\lb\lambda_{N_1+1},...,\lambda_{N_2}\rb/I_2
\end{array}$$
we will follow Definition 5.4.2 of \cite{RoSPS} and define $B_1\otimes_A^sB_2$ the {\em separated tensor product of $B_1$ and $B_2$ over $A$} to be the ring 
$$B_1\otimes_A^s B_2:=A\ls y_1,...,y_{M_2}\rs\lb\lambda_1,...,\lambda_{N_2}\rb/(I_1\cup I_2).$$
By Theorem 5.2.6 of \cite{RoSPS}, this product is independent of the presentations of $B_1$ and $B_2$.

Now following the Definition 5.3.1 of \cite{RoSPS}:

\begin{defn}\label{grf}
 A {\em generalized ring of fractions} over $\smn$ is inductively defined as follows:

i) $\smn$ is a generalized ring of fractions over $\smn$.

ii) If $B'$ is a generalized ring of fractions over $\smn$, and  $g,\:f_1,...,f_s,f'_1,...,f'_t\in B'$, then 
$$B=B\ls f_1/g,...,f_s/g\rs\lb f'_1/g,...,f'_t/g\rb$$
is also a generalized ring of fractions over $\smn$.
\end{defn}

If $B$ is a generalized ring of fractions over $\smn$ then the members of $B$ can be treated as functions over the domain of $B$ which is defined below.

\begin{defn}\label{dom}
Let $K'$ be a complete extension of $K$ and let 
$$B=S_{m+s,n+t}/(\{g_ix_i-f_i\}_{i=m+1}^{m+s}\cup\{g'_j\rho_j-f'_j\}_{j=n+1}^{n+t})$$ 
be a generalized ring of fractions over $\smn$. Let $\bar{g}_i$, $\bar{f}_i$, $\bar{g}'_j$ and $\bar{f}'_j$ denote the images of $g_i$, $f_i$, $g'_j$ and $f'_j$ in $B$ respectively, then we define $\kdom B$, the {\em $K'$-rational points in the domain of $B$} to be the set
$$\{\pb\in(K^\circ)^m\times(K^{\circ\circ})^n: |\bar{f}_i(\pb)|\leq|\bar{g}_i(\pb)|\neq 0, |\bar{f}'_j(\pb)|<|\bar{g}'_j(\pb)|\text{ for all }i,j\}.$$
When $K'=K$ we will simply write $\dom B$ instead of $K\text{-Dom}_{m,n}B$.
\end{defn}
Note that for $\bar{K}$ the completion of algebraic closure of $K$, our definition for $\bar{K}$-rational points of the domain of $B$ coincides with what is called the {\em domain of $B$} in \cite{modcomp}.

For a generalized ring of fractions $B$ over $\smn$ as above and an ideal $I\subset B$ we will use the customary notation $\dom B \cap V(I)_K$ to denote the set
$$\dom B \cap V(I)_K:=\{\pb\in\dom B: f(\pb)=0\text{ for all }f\in I\},$$
and omit the subscript $K$ when $K$ is algebraically closed. 

With these notations established, we can give an alternate definition of a $D$-semianalytic set which we find useful because it connects algebraic objects with geometric objects in a customary way.
\begin{defn}\label{Dsemi}
A $D$-semianalytic subset of $(K^\circ)^{m+n}$ is a finite union of sets of the form $\dom B\cap V(I)_K$ where $B$ is a generalized ring of fractions over $\smn$ and $I$ is an ideal of $B$.
\end{defn}

Next, we will make observations which will be helpful in the proofs of Theorems \ref{parass} and \ref{densesub}. Let $B$ be as in Definition \ref{dom} and let $J\subset S_{m+s,n+t}$ be the ideal $(\{g_ix_i-f_i\}_{i=m+1}^{m+s}\cup\{g'_j\rho_j-f'_j\}_{j=n+1}^{n+t})$. With this ideal we can associate a $K$-rational variety $V(J)_K$ of $(K^\circ)^{m+s}\times(K^{\circ\circ})^{n+t}$ and a semianalytic set 
$$X:=\{\pb\in (K^\circ)^{m+s}\times(K^{\circ\circ})^{n+t}:\pb\in V(J)_K\text{ and }\prod_{i=m+1}^{m+s} \bar{g}_i(\pb)\cdot\prod_{j=n+1}^{n+t} \bar{g}'_j(\pb)\neq 0\}.$$ 
Then the coordinate projection $\pi:(K^\circ)^{m+s}\times(K^{\circ\circ})^{n+t}\rightarrow (K^\circ)^m\times(K^{\circ\circ})^n$ maps $X$ bijectively onto $\dom B$. Furthermore by the Implicit Function Theorem, in fact if $X$ is non-empty then it is an $(m+n)$-dimensional $K$-analytic manifold (see Definition \ref{mani}). 

In order to understand properties of the set $\dom B$ we will often look at the set $X$. However to avoid complications coming from working over a not necessarily algebraically closed field and to avoid irreducible components of $X$ that are contained in the excluded set 
$$Z:=\{\pb\in(K^\circ)^{m+s}\times(K^{\circ\circ})^{n+t}: \prod_{i=m+1}^{m+s} g_i(\pb)\cdot \prod_{j=n+1}^{n+t} g'_j(\pb) =0\}$$ 
we will prefer to work with the largest ideal that vanishes on $\dom B \cap V(I)_K$. That is, the ideal 
$$\ideal(\dom B \cap V(I)_K):=\{f\in B: f(\pb)=0\text{ for all }\pb\in\dom B \cap V(I)_K\}$$
whose corresponding ideal $J'$ in $S_{m+s,n+t}$ does not have a minimal prime divisor $\mfp$ such that $V(\mfp)_K\subset Z$, will show up often in our arguments.

There is a special type of $D$-semianalytic sets that will attract our attention.

\begin{defn}\label{Rdom}
Let $B$ be a generalized ring of fractions over $\smn$ such that at each inductive step of the construction of $B$ as described in Definition \ref{grf}, the ideal $(g,f_1,...,f_s,f'_1,...,f'_t)$ is the unit ideal of $B'$. Let $\bar{K}$ be an algebraically closed complete extension of $K$, then $\bar{K}\text{-Dom}_{m,n}B$ is called an {\em $R$-subdomain of $(\bar{K}^\circ)^{m+n}$}.
\end{defn}

$R$-domains generalize the rational domains of the affinoid geometry (as given in Definition 5 of $\S$7.2.3 of \cite{bgr}). They are also examples of quasi-affinoid subdomains and have the universal property described in Definition 5.3.4 of \cite{RoSPS}. This implies that if $A$ and $B$ are two generalized rings of fractions over $\smn$ such that $X:=\bar{K}\text{-Dom}_{m,n}A=\bar{K}\text{-Dom}_{m,n}B$ is an $R$-domain, then $A\simeq B$. This enables us to define a ring of analytic functions on $X$ with coefficients from $K$, which we will denote by $\mathcal{O}(X)_K:=A\simeq B$.

\section{Normalization and Dimension Theory}

In this section we briefly discuss the dimension theory of $D$-semianalytic sets and Parameterized Normalization Lemma as they are some of our main tools in the proofs of the main results of this paper. These were studied extensively in \cite{dimtheory} and we will start by the basic definitions and results from that source. However we will slightly improve on these results to put them in a form that enables us to prove the main theorems of this paper.

For an arbitrary subset of $K^m$ we can define several notions of dimension as follows.

\begin{defn}\label{gdim}
 Define the {\em geometric dimension}, g-dim $X$,  of a nonempty set $X\subset K^m$ to be the greatest integer $d$ such that the image of $X$ under coordinate projection onto a $d$-dimensional coordinate hyper-plane has an interior point. For $X=\emptyset$ define g-dim $X=-1$. 

The {\em weak dimension}, w-dim $X$, of a nonempty set $X\subset K^m$ to be the greatest integer $d$ such that the image of $X$ under coordinate projection onto a $d$-dimensional coordinate hyper-plane is somewhere dense. For $X=\emptyset$ define g-dim $X=-1$. 

For a generalized ring of fractions $B$ and an ideal $I$ of $B$, k-dim $B/I$ will denote the Krull dimension of the algebra $B/I$.
\end{defn}

Note that the reference to a $d$-dimensional coordinate hyperplane in the above definition is superfluous as one may as well work with any $d$-dimensional hyperplane. One of the main results of \cite{dimtheory} was

\begin{thm}\label{alleq0}
Assume {\em Char }$K=0$. For a generalized ring of fractions $B$ over $\smn$ and ideal $I$ of $B$ satisfying $I=\mathcal{I}(\dom B\cap V(I)_K)$,
$$\text{w-dim }\dom B\cap V(I)_K=\text{g-dim }\dom B\cap V(I)_K=\text{ k-dim }B/I.$$
\end{thm}

Let $B=\Smn/J$ be a generalized ring of fractions over $\smn$, $I$ be an ideal of $B$, and $K'$ be a complete extension of $K$. It is natural to ask the relation between the geometric dimensions of $\dom B\cap V(I)_K$ and $\kdom B'\cap V(I')_{K'}$ where $B':=S_{0,0}(E,K')\otimes^s_{S_{0,0}(E,K)}B=\Smn(E,K')/J\cdot\Smn(E,K')$ and $I'=I\cdot B'$. In other words, if $B'$ is obtained from $B$ through the ground field extension from $(E,K)$ to $(E,K')$, is it true that if $\varphi(x,\rho)$ is a quantifier-free formula of the language of \cite{rigid} or \cite{modcomp} then the $D$-semianalytic set $X=\{\pb\in (K^\circ)^m\times (K^{\circ\circ})^n:\varphi(\pb)\}$ and the $D$-semianalytic set $X'=\{\pb\in ({K'}^\circ)^m\times ({K'}^{\circ\circ})^n:\varphi(\pb)\}$ have the same geometric dimension. If we avoid the obvious complications that may arise from the fact that $K$ is not necessarily algebraically closed, the answer is ``yes''. 

\begin{thm}\label{extensiondim}
Let {\em Char $K=0$}, $B$, $I$, $J$, $K'$ and $B'$ be as above, and assume that $I=\ideal(\dom B\cap V(I)_K)$, then
$$\text{g-dim }\dom B\cap V(I)_K=\text{g-dim }\kdom B'\cap V(I')_{K'}.$$
\end{thm}
\begin{proof}
Let $I^*=\ideal(\kdom B'\cap V(I')_{K'})$ and notice that by Theorem \ref{alleq0} it is enough to prove that k-dim $B/I=$ k-dim $B'/I^*$. Let $\bar{I}$ be the ideal that corresponds to $I$ in $\smn(E,K)$ and let $\bar{I}^*$ be the ideal that corresponds to $I^*$ in $\smn(E,K')$ so that $\bar{I}\subset\bar{I}^*\cap\smn(E,K)$. Let us write 
$$B=\Smn/(\{g_ix_i-f_i\}_{i=m+1}^{m+M}\cup\{g'_j\rho_j-f'_j\}_{j=n+1}^{n+N})$$
so that for any point $\pb\in \dom B\cap V(I)_K$ there corresponds a maximal ideal $\mfm$ of $\smn(E,K)$ such that $g_i,g'_j\not\in\mfm$ for any $i,j$. Note that by Proposition 5.4.10 of \cite{RoSPS} $\smn(E,K')$ is faithfully flat over $\smn(E,K)$ and therefore for any such $\mfm$ there exists a maximal ideal $\mfm'$ of $\smn(E,K')$ such that $g_i,g'_j\not\in\mfm'$ for any $i,j$ and $\mfm'\cap\smn(E,K)=\mfm$. By construction of $\bar{I}^*$ we have $\mfm'\supset \bar{I}^*$ and hence $\mfm\supset \bar{I}^*\cap\smn(E,K)$ for any such maximal ideal $\mfm$. Now by the Nulstellensatz (Theorem 4.1.1 of \cite{RoSPS}) we have $\bar{I}^*\cap\smn(E,K)=\bar{I}$. On the other hand, once again by the faithful flatness of $\smn(E,K')$ over $\smn(E,K)$ we have k-dim $B/I=$ k-dim $\smn(E,K)/\bar{I}=$ k-dim $\smn(E,K')/(\bar{I}^*\cap\smn(E,K'))=$ k-dim $B'/I^*$ and the result follows.
\end{proof}

The above theorem allows us to prove another result in dimension theory of $D$-semianalytic sets.

\begin{thm}\label{bndry}
Let {\em Char $K=0$} and $X\subset(K^\circ)^m\times(K^{\circ\circ})^n$ be $D$-semianalytic. Let us write $\bar{X}$ for the closure of $X$ in the metric topology, then there is a $D$-semianalytic set $Y$ such that $\bar{X}\setminus X\subset Y$ and g-dim $Y<$ g-dim $X$.
\end{thm}
\begin{proof}
By Theorem \ref{alleq0} it is enough to consider the case $X=\dom B\cap V(I)_K$ where $B$ is a generalized ring of fractions over $\smn$ and $I\subset B$ is a prime ideal satisfying $I=\ideal(\dom B\cap V(I)_K)$. First we will prove the theorem in the case that $K$ is algebraically closed. We would like to point out that this case was already proved in Theorem 4.3 of \cite{strat} by a somewhat different argument. Note that in this case, by the Quantifier Elimination Theorem of \cite{rigid}, one can choose $Y=\bar{X}\setminus X$ as it is $D$-semianalytic. 

Write 
$$B=\Smn/(\{g_ix_i-f_i\}_{i=m+1}^{m+M}\cup\{g'_j\rho_j-f'_j\}_{j=n+1}^{n+N})=\Smn/J.$$
Let  $\pi:(K^\circ)^{m+M}\times(K^{\circ\circ})^{n+N}\rightarrow (K^\circ)^m\times(K^{\circ\circ})^n$ be the coordinate projection, $\bar{I}$ be the ideal of $\Smn$ that corresponds to $I$ and let $X'\subset (K^\circ)^{m+M}\times(K^{\circ\circ})^{n+N}$ be the semianalytic set $V(\bar{I})_K\cap\{\qb:\prod_{i=m+1}^{m+M} g_i(\qb)\cdot\prod_{j=n+1}^{n+N}g_j'(\qb)\neq 0\}$, so that $\pi(X')=X$. Note that in this case g-dim $X=$ k-dim $\Smn/\bar{I}$.

Now let $\pb\in\bar{X}\setminus X$ and assume that $\pb\in\bar{X}\cap\dom B$, then there is a $\qb\in X'$ such that $\pi(\qb)=\pb$ contradicting $\pb\not\in X$. Therefore $(\bar{X}\setminus X)\cap \dom B=\emptyset$. On the other hand by Lemma 6.3 of \cite{uniform}, there is a $\qb\in V(\bar{I})_K$ such that $\pi(\qb)=\pb$ and hence it must be the case that $\qb\in V(I\cup\{\prod_{i=m+1}^{m+M} g_i\cdot\prod_{j=n+1}^{n+N}g_j'\})_K$ and g-dim $V(I\cup\{\prod_{i=m+1}^{m+M} g_i\cdot\prod_{j=n+1}^{n+N}g_j'\})_K<$ g-dim $X$. Setting 
$$Y:=\pi(V(I\cup\{\prod_{i=m+1}^{m+M} g_i\cdot\prod_{j=n+1}^{n+N}g_j'\})_K)$$
finishes the case where $K$ is algebraically closed.

For an arbitrary complete field $K$, let $\bar{K}$ be an algebraically closed complete extension of $K$ and let $\bar{B}$ be the generalized ring of fractions $\Smn(E,\bar{K})/(J\cdot\Smn(E,\bar{K}))$. Put 
$$X'':=\bar{K}\text{-Dom}_{m,n}\bar{B}\cap V(I\cdot\bar{B})_{\bar{K}}$$
and notice that by Theorem \ref{extensiondim} g-dim $X=$ g-dim $X''$. Notice also that $\bar{X}\setminus X$ is contained in the set $Y:=(\bar{X}''\setminus X'')\cap (K^\circ)^{m}\times(K^{\circ\circ})^{n}$. Observe that by the Quantifier Elimination Theorem (Corollary 4.3) of \cite{modcomp}, the set $(\bar{X}''\setminus X'')$ is $D$-semianalytic which can be written in the form
$$\bar{X}''\setminus X''=\bigcup_i\bar{K}\text{-Dom}_{m,n}B_i\cap V(I_i)_{\bar{K}}$$
where each $B_i$ is a generalized ring of fractions over $\smn(E,K)$ and $I_i\subset B_i$ is an ideal. This shows that the set $Y$ which is also the same set as $\bigcup_i\dom B_i\cap V(I_i)_K$ is in fact $D$-semianalytic. Now the result follows from Theorem \ref{extensiondim} as 
$$\text{g-dim }Y\leq\text{ g-dim }\bar{X}''\setminus X''<\text{ g-dim }X''=\text{ g-dim }X.$$
\end{proof}

One of the main tools for proving Theorem \ref{alleq0} (Theorem 6.2 of \cite{dimtheory}) as well as obtaining the results in this paper is the Parameterized Normalization Lemma (Lemma 5.3 of \cite{dimtheory}). It is known that given a quasi-affinoid algebra $B$, it is not always possible to find a ring of separated power series $\smn$ such that there is a finite injection $\phi:\smn\rightarrow A$ (see Example 2.3.5 of \cite{RoSPS}). Nevertheless we can break up the $D$-semianalytic set associated with $B$ into finitely many smaller $D$-semianalytic sets whose associated quasi-affinoid algebras can be normalized in the sense above. Furthermore in this process of breaking-up, several key properties are preserved including the parameter structure of those algebras. Note that given a quasi-affinoid algebra there may be more than one way of considering it as a generalized ring of fractions depending on which variables we choose to represent the coordinates of the space in which the $D$-semianalytic set lives and which variables correspond to fractions, although this fact plays no role in our discussions. 

Let us write $\Smn$ for the ring of separated power series in the variables $(x_1,...,x_m)$, $(\rho_1,...,\rho_n)$, $(y_1,...,y_M)$ and $(\lambda_1,...,\lambda_N)$. Now suppose $B$ is a generalized ring of fractions over $\Smn$, then assigning appropriate $x$ or $\rho$ variable to fractions involving terms built up inductively from only the variables $x_1,...,x_m$ and $\rho_1,...,\rho_n$, $B$ can be written in the form
\begin{multline*}\label{paragrf}\tag{$\ast$}
B=(\smn\ls x_{m+1},....,x_{m+s}\rs\lb\rho_{n+1},...,\rho_{n+t}\rb/J_1)\\
\ls y_{1},...,y_{M+S}\rs\lb\lambda_{1},...,\lambda_{N+T}\rb/J_2,
\end{multline*}
We will call the $x$ and $\rho$ variables appearing above the {\em parameter variables}. 

For a ring of separated power series $S_{m+M,n+N}$, and a fixed $1\leq j\leq M$, we will call an automorphism $\phi$ of $S_{m+M,n+N}$ which is defined by $\phi(y_i):=y_i+y_j^{r_i}$ for some $r_i\in\mathbb{N}^+$ for $i<j$ and $\phi(y_i)=y_i$ for $i\geq j$, a Weierstrass change of variables among $y$ variables. We define the Weierstrass change of variables among $x$, $\rho$ or $\lambda$ variables similarly. We will use the same term for compositions of such variable changes also, so that a Weierstrass change of variables among a single type of variables ($x$, $\rho$, $y$ or $\lambda$) respects the sort of the variables as well as the parameter structure.

In the following discussions, we will often come across a situation where for a quasi-affinoid algebra $B=\Smn/J$, and an ideal $I\subset B$, the corresponding ideal $\bar{I}$ of $\Smn$ is such that after a  Weierstrass variable change $\phi$, we can find another quasi affinoid algebra $A$ such that $A\subset\Smn/\phi(\bar{I})$ is a finite inclusion. In such cases we will abuse the notation and say $A\subset B/\phi(I)$ is a finite inclusion.

Combining Lemmas 5.3 and 5.5 of \cite{dimtheory} we have:

\begin{lem}[Parameterized Normalization Lemma]\label{parnorm1}
Let $B$ be as in \eqref{paragrf} and $I$ be an ideal of $B$, then there exist finitely many generalized rings of fractions 
\begin{multline*}
B_i=(\smn\ls x_{m+1},....,x_{m+s_i}\rs\lb\rho_{n+1},...,\rho_{n+t_i}\rb/J_{i,1})\\
\ls y_{1},...,y_{M+S_i}\rs\lb\lambda_{1},...,\lambda_{N+T_i}\rb/J_{i,2},
\end{multline*}
with parameter rings $A_i:=\smn\ls x_{m+1},....,x_{m+s_i}\rs\lb\rho_{n+1},...,\rho_{n+t_i}\rb/J_{i,1}$, ideals $I_i\subset B_i$, and integers $M_i$, $N_i$ such that

i) $V(I)_K\cap\Dom B=\bigcup_i(V(I_i)\cap\Dom B_i)$,

ii) $I_i=\ideal(V(I_i)_K\cap\Dom B_i)$, k-dim $B/I\geq$ k-dim $B_i/I_i$,

iii) after a Weierstrass change of variables $\phi_i$ among $y$ and $\lambda$ variables separately we have a finite inclusion
$$A_i/(A_i\cap I_i)\ls y_1,...,y_{M_i}\rs\lb\lambda_1,...,\lambda_{N_i}\rb\subset B_i/\phi_i(I_i).$$
\end{lem}

Although this statement of Normalization is quite useful for working with the geometric properties of projections of $D$-semianalytic sets, of course one can carry the process of normalization one more step to make it look more like other well known normalization results from algebra.

\begin{lem}\label{morenorm} 
Let $B$ be as in \eqref{paragrf} and $I$ be an ideal of $B$, then there exist finitely many generalized rings of fractions 
$$B_i=S_{m+s_i+M+S_i,n+t_i+N+T_i}/J_i,$$ 
in the variables $(x_1,...,x_{m+s_i})$, $(\rho_1,...,\rho_{n+t_i})$, $(y_1,...,y_{M+S_i})$, $(\lambda_1,...,\lambda_{N+T_i})$; ideals $I_i$, and integers $m_i$, $n_i$, $M_i$, $N_i$ such that 

i) $V(I)_K\cap\Dom B=\bigcup_i(V(I_i)_K\cap\Dom B_i)$,

ii) $I_i=\ideal(V(I_i)_K\cap\Dom B_i)$, 

iii) $M_i+N_i\leq M+N$, $m_i+n_i\leq m+n$ and after a Weierstrass change of variables $\psi_i$ among $x$, $\rho$, $y$ and $\lambda$ variables separately we have a finite inclusion
$$S_{m_i+M_i,n_i+N_i}\rightarrow B_i/\psi_i(I_i).$$
\end{lem}

\begin{proof}
After applying Lemma \ref{parnorm1} we may assume that there are integers $M'$ and $N'$, a generalized ring of fractions $A$ over $\smn$ and  a Weierstrass change of variables $\phi$ among $y$ and $\lambda$ variables separately such that 
$$A/(I\cap A)\ls y_1,...,y_{M'}\rs\lb\lambda_1,...,\lambda_{N'}\rb\rightarrow B/\phi(I)$$
is a finite inclusion. Apply Lemma \ref{parnorm1} once more to $A/(I\cap A)$ to get generalized rings of fractions $A_{i}$ over $\smn$, ideals $J_{i}'\subset A_{i}$, integers $m_{i}$, $n_{i}$ and Weierstrass changes of variables $\phi_{i}$ among $x$ and $\rho$ variables separately such that 
$$V(I\cap A)_K\cap\dom A=\bigcup_j(\dom A_{i}\cap V(J_{i}')_K),$$
$$J_{i}'=\ideal( \dom A_{i}\cap V(J_{i}')_K),\text{ k-dim }A/(I\cap A)\geq \text{ k-dim }A_{i}/J_{i}',$$
 for all $i$ and $S_{m_{i},n_{i}}\subset A_{i}/\phi_i(J_{i}')$ is a finite inclusion. Then by Lemmas 2.4  and 2.5 of \cite{dimtheory} we have 
$$A_{i}/J_{i}'\ls y_1,...,y_{M'}\rs\lb\lambda_1,...,\lambda_{N'}\rb\subset B_{i}/\phi(I\cup J_{i}'\cdot B_{i}),$$ 
is a finite inclusion where $B_{i}$ is the generalized ring of fractions over $\Smn$ which is the separated tensor product $A_{i}\otimes^s_{A}B$. Again by the same lemma
$$S_{m_{i}+M',n_{i}+N'}\subset B_{i}/\phi_{i}\circ\phi(I\cup J_{i}'\cdot B_{i}).$$
is also a finite inclusion.

Let $I_{i}=\ideal(V((I\cup J_{i}')\cdot B_{i})_K\cap\Dom B_{i})$, and notice that in case k-dim $B_{i}/I_{i}=m_{i}+n_{i}+M'+N'$ then we can replace the ideal $(I\cup J_{i}')B_{i}$ with $I_{i}$ and still have the above map $\phi_{i}\circ\phi$ a finite inclusion. Otherwise, the Krull dimension goes down and we proceed inductively applying the process described above to the generalized ring of fractions $B_{i}$ and ideal $I_{i}$. Note that by Lemma 5.5 of \cite{dimtheory} we have $m_{i}+n_{i}\leq m+n$ and $M'+N'\leq M+N$.
\end{proof}

\begin{rem}\label{respardense}

i) Let us write $B_i=S_{m+s_i+M+S_i,n+t_i+N+T_i}/J_i$ for the generalized rings of fractions we found in the previous lemma. One important property in the above setting is that $x$, $\rho$, $y$ and $\lambda$ variables are not mixed under $\phi_{i}\circ\phi$, i.e. the image under the change of variables (automorphisms) $\phi_{i}\circ\phi$ of $S_{m+s_i,n+t_i}$ is again $S_{m+s_i,n+t_i}$. We will call such an automorphism a {\em parameter respecting automorphism} of quasi-affinoid algebras with a parameter structure. 

ii) Another observation we are going to make about the parameterized normalization will be useful in the proof of Theorem \ref{densesub}. Suppose that the projection of $\Dom B\cap V(I)_K$ onto the coordinate hyperplane $(K^\circ)^m\times(K^{\circ\circ})^n$ is somewhere dense and observe that this projection is contained in $\bigcup_i(\dom A_{i}\cap V(J_{i}')_K)$ where each $A_{i}$ is a generalized ring of fractions over $\smn$ and each $J_{i}'\subset A_{i}$ is an ideal as in the proof above. Then for some index $i_0$ the projection of the $D$-semianalytic subset $\Dom B_{i_0}\cap V(I_{i_0})_K$ onto $(K^\circ)^m\times(K^{\circ\circ})^n$ which is contained in $\dom A_{i_0}\cap V(J_{i_0}')_K$ contains a somewhere dense subset and hence by Lemma 4.9 of \cite{dimtheory} we have k-dim $A_{i_0}/J_{i_0}'=m+n$, and therefore $m_{i_0}+n_{i_0}=m+n$.
\end{rem}

We will finish this section with one last observation about a ``localization'' property of normalization.

\begin{lem}\label{local}
Let $B$ be a generalized ring of fractions over $\smn$ and $I$ be an ideal of $B$. Suppose that $S_{m',n'}\subset B/I$ is a finite inclusion, the origin $\mathbf{0}$ is in  $\dom B\cap V(I)_K$ and the maximal ideal $\mfm$ of $B$ corresponding to the origin $\mathbf{0}$ is such that k-dim $B_{\mfm}/IB_{\mfm}=m'+n'$. Then for all $\varepsilon\in K^{\circ}$, $\varepsilon\neq 0$ such that the open ball $B_K(\mathbf{0},|\varepsilon|)$ of center $\mathbf{0}$, radius $|\varepsilon|$ is contained in $\dom B$, there is a $\delta\in K^{\circ}$, $\delta\neq 0$ such that 
$$T_{m'+n'}(\delta)\subset T_{m'+n'}(\delta)\ls x_{m'+1}/\varepsilon,...,x_m/\varepsilon,\rho_{n'+1}/\varepsilon,...,\rho_n/\varepsilon\rs/I $$
is a finite inclusion, where $T_{m'+n'}(\delta)$ stands for the Tate ring 
$$K\ls x_1/\delta,...,x_{m'}/\delta,\rho_1/\delta,...,\rho_{n'}/\delta\rs.$$
\end{lem}
\begin{proof}
For simplicity in notation we will assume $n'=n=0$ but the arguments below will also work in the general case. Let 
$$f_{m'+i}=x_{m'+i}^{s_{m'+i}}+b_{m'+i,s_{m'+i}-1}(x_1,...,x_{m'})x_{m'+i}^{s_{m'+i}-1}+...+b_{m'+i,0}(x_1,...,x_{m'})$$ 
be the lowest degree monic polynomial in $I\cap K\ls x_1,...,x_{m'}\rs\left[x_{m'+i}\right]$ for $i=1,..,m-m'$. For each such $i$ and $j=1,...,s_{m'+i}$ write
$$b_{m'+i,s_{m'+i}-j}=c_{m'+i,s_{m'+i}-j}+d_{m'+i,s_{m'+i}-j},$$
where $=c_{m'+i,s_{m'+i}-j}\in K$ and $d_{m'+i,s_{m'+i}-j}\in (x_1,...,x_{m'})K\ls x_1,...,x_{m'}\rs$. By continuity there is a $\delta\in K^{\circ}$, $\delta\neq 0$ such that whenever $\pb\in B_K(\mathbf{0},|\delta|)\cap(K^\circ)^m$ we have $|d_{m'+i,s_{m'+i}-j}(\pb)|<|\varepsilon|^j$ for all $i=1,...,m-m'$ and $j=1,...,s_{m'+i}$.

Notice that in this case each $f_{m'+i}$ is regular in $x_{m'+i}/\varepsilon$ of some degree $r_{m'+i}$, $1\leq r_{m'+i}\leq s_{m'+i}$, in 
$$B\ls x_1/\delta,...,x_{m'}/\delta,x_{m'+1}/\varepsilon,...,x_m/\varepsilon\rs\simeq T_{m'}(\delta)\ls x_{m'+1}/\varepsilon,...,x_m/\varepsilon\rs.$$ 
This shows that $T_{m'}(\delta)\ls x_{m'+1}/\varepsilon,...,x_m/\varepsilon\rs/I$ is finite over $T_{m'}(\delta)$. 
 On the other hand by Lemma 3.7 of \cite{dimtheory} k-dim $T_{m'}(\delta)\ls x_{m'+1}/\varepsilon,...,x_m/\varepsilon\rs/I=m'$ and we see that
$$T_{m'}(\delta)\hookrightarrow T_{m'}(\delta)\ls x_{m'+1}/\varepsilon,...,x_m/\varepsilon\rs/I$$
is a finite inclusion.
\end{proof}

\section{Parameterized Stratification}

In this section we will sharpen the Smooth Stratification Theorem of \cite{dimtheory} and \cite{strat} for $D$-semianalytic sets using the normalization results of the previous section. Our goal is to prove that a $D$-semianalytic set is in fact a finite union of $D$-semianalytic manifolds which remain manifolds when specialized at a given point of the parameter space. The next lemma handles the main step in proving the above statement by establishing that a properly normalized quasi-affinoid variety can be written as a union of an analytic manifold and a smaller dimensional variety where the manifold is locally the graph of some analytic functions, the normalized parameters being functions of only the free parameters of the normalization.

\begin{lem}\label{parajac}
Let {\em Char $K=0$} and let $I$ be a prime ideal of $\Smn$ such that $I=\ideal(V(I)_K)$. Suppose $S_{m'+M',n+N'}\subset \Smn/I$ is a finite inclusion respecting parameters (see Remark \ref{respardense} (i)) then we have
$$V(I)_K=Z\cup Y,$$
where $Z$ is a $(m'+M'+n'+N')$-dimensional analytic manifold whose charts are given by projections onto the space $(K^\circ)^{m'+M'}\times(K^{\circ\circ})^{n'+N'}$ and where $Y=V(I')_K$ for some ideal $I'\supset I$.

More precisely, for each $\pb\in Z$ there is an open neighborhood $U$ of $\pb$ and analytic functions $\alpha_{m'+i}(x',\rho')$, $\beta_{n'+j}(x',\rho')$, $\gamma_{M'+k}(x',\rho',y',\lambda')$, $\delta_{N'+l}(x',\rho',y',\lambda')$ of $U$ for $i=1,..,m-m'$, $j=1,...,n-n'$, $k=1,...,M-M'$, $l=1,...,N-N'$ and multi variables $x'=(x_1,...,x_{m'})$, $\rho'=(\rho_1,...,\rho_{n'})$, $y'=(y_1,...,y_{M'})$, $\lambda'=(\lambda_1,...,\lambda_{N'})$ such that
\begin{multline*}
Z\cap U= V(\{x_{m'+i}-\alpha_{m'+i}\}_{i=1}^{m-m'}\cup\{\rho_{n'+j}-\beta_{n'+j}\}_{j=1}^{n-n'}\cup\\
\{y_{M'+k}-\gamma_{M'+k}\}_{k=1}^{M-M'}\cup\{\lambda_{N'+l}-\delta_{N'+l}\}_{l=1}^{N-N'})_K\cap U.
\end{multline*}
\end{lem}

\begin{proof}
 For simplicity in notation we will assume that $N=N'=n=n'=0$, nevertheless our arguments will be valid in the general case. 

 Let $p_{m'+i}$ be the unique lowest degree monic polynomial in $I\cap K\ls x\rs\left[x_{m'+i}\right]$ for each $i$ and $q_{M+j}$ be the lowest degree monic polynomial in $I\cap K\ls x',y'\rs\left[y_{M'+j}\right]$ for each $j$. Define
\begin{displaymath}\Delta:=\text{det}\left(\displaystyle\frac{\partial \{p_{m'+i},q_{M'+j}\}_{i,j}}{\partial \{x_{m'+i},y_{M'+j}\}_{i,j}}\right)=\displaystyle\frac{\partial p_{m'+1}}{\partial x_{m'+1}}\cdots\frac{\partial p_{m}}{\partial x_{m}}\frac{\partial q_{M'+1}}{\partial y_{M'+1}}\cdots\frac{\partial q_{M}}{\partial y_{M}},
\end{displaymath}
as the Jacobian matrix above is upper triangular. Define also 
$$Z:=\{\pb\in V(I)_K:\Delta(\pb)\neq 0\}.$$

Notice that because each $p_{m'+i}$ is the lowest degree monic polynomial and $I$ is a prime ideal, we have $\Delta\not\in I$. Therefore by Theorem 6.2 of \cite{dimtheory} the geometric dimension of $Y:=V(I)_K\setminus Z=V(I\cup\{\Delta\})_K$ is less than $m'+M'+n'+N'$.

On the other hand, the fact that for all $\pb\in Z$, we have $\partial p_{m'+i}/\partial x_{m'+i}(\pb)\neq 0$, $\partial q_{M'+j}/\partial y_{M'+j}(\pb)\neq 0$ for all $i,j$ implies that each $p_{m'+i}$ is regular of degree one in $x_{m'+i}$ and each $q_{M'+j}$ is regular of degree one in $y_{M'+j}$ in a rational neighborhood $W$ of $\pb$. Here the regularity of these polynomials is in the sense of Definition 2.3.7 of \cite{RoSPS}. Furthermore $p_{m'+i}$, $q_{M'+j}$ generate $I\mathcal{O}(U)$ for some neighborhood $U\subset W$ of $\pb$ by Theorem 30.4 of \cite{matsu}. Now the statement follows from the Weierstrass Division Theorem (Theorem 2.3.8 of \cite{RoSPS}).
\end{proof}

In order to prove Parameterized Smooth Stratification Theorem, first we need to clarify what we mean by a $D$-semianalytic manifold. Following \cite{rigid} and \cite{serre}, the definition is as follows.

\begin{defn}\label{mani}
A $D$-semianalytic subset $X$ of $K^m$ is a $d$-dimensional $K$-analytic manifold if it is endowed with a system of charts $(U,\varphi_U:U\rightarrow K^d)$ such that the transition maps $\varphi_U\circ\varphi_V^{-1}$ are locally given by convergent power series over $K$. Note that here the topology of $X$ is the subspace topology inherited from the metric topology of $K$.
\end{defn}

In what follows, we will be interested in specializations of $D$-semianalytic subsets of $(K^{\circ})^{m+M}\times(K^{\circ\circ})^{n+N}$ at points from $(K^{\circ})^{m}\times(K^{\circ\circ})^{n}$ which we consider as the parameter space. The Parameterized Smooth Stratification Theorem for $D$-semianalytic sets is the following:

\begin{thm}[Parameterized Smooth Stratification]\label{parass}
Let {\em Char $K=0$} and let $X$ be a $D$-semianalytic subset of $(K^{\circ})^{m+M}\times(K^{\circ\circ})^{n+N}$, then we can write $X$ as a finite union of $D$-semianalytic manifolds $X_i$ such that 
$$\sup_i\text{ g-dim }X_i=\text{ g-dim }X$$
and for each $\pb\in(K^{\circ})^{m}\times(K^{\circ\circ})^{n}$, the fiber $X_i(\pb)$ of $X_i$ at $\pb$ is either empty or is a $D$-semianalytic manifold.
\end{thm}
\begin{proof}
We may assume that $X=\Dom B\cap V(I)_K$ for some generalized ring of fractions $B$ over $\Smn$ and ideal $I$ of $B$. By Lemma \ref{morenorm} we may also assume that there are integers $m',n',M',N'$ and a Weierstrass automorphism $\phi$ respecting the parameters (see Remark \ref{respardense}(i)) such that $S_{m'+M',n'+N'}\subset B/\phi(I)$ is a finite inclusion and  $I=\ideal(\Dom B\cap V(I)_K)$. By Lemma 4.5 of \cite{dimtheory} each minimal prime divisor $\mfp$ of $I$ satisfies $\mfp=\ideal(\text{Dom}_{m+M}B\cap V(\mfp)_K)$ so we may also assume $I$ is prime.

The proof will be by induction on g-dim $X:=d=m'+M'+n'+N'$. There is nothing to prove if $X=\emptyset$ or if $d=0$, as in that case we have the statement by Remark 3.2 and Theorem 6.2 of \cite{dimtheory}. Therefore we will assume that $d>0$. 

The following maps between the affinoid spaces will play an important role in understanding the relations between geometric objects we are concerned about. Write $B=S_{m+s+M+S,n+t+N+T}/J$ and notice that $\phi$ induces a one-to-one analytic transformation 
$$\phi':(K^\circ)^{m+s+M+S}\times(K^{\circ\circ})^{n+t+N+T}\rightarrow (K^\circ)^{m+s+M+S}\times(K^{\circ\circ})^{n+t+N+T}$$
with non-zero Jacobian so that $\phi'(V(\bar{I})_K)=V(\phi(\bar{I}))_K$ where $\bar{I}$ is the ideal of $S_{m+s+M+S,n+t+N+T}$ that corresponds to $I$. Another map related to our construction is the projection map 
$$\pi:(K^\circ)^{m+s+M+S}\times(K^{\circ\circ})^{n+t+N+T}\rightarrow (K^\circ)^{m+M}\times(K^{\circ\circ})^{n+N}$$
which is one-to-one when restricted to $\pi^{-1}(\Dom B)\cap V(J)_K$.

Note that in our setting there is an open set $U\subset (K^\circ)^{m+s+M+S}\times(K^{\circ\circ})^{n+t+N+T}$ such that $X=\pi(V(\bar{I})_K\cap U)$ and $\bar{I}=\ideal(V(\bar{I})_K)$. On the other hand by Theorem \ref{alleq0} we have g-dim $V(\bar{I})_K\cap U=d$. Next, apply Lemma \ref{parajac} to $\phi'(V(\bar{I})_K)=V(\phi(\bar{I}))_K$ (replacing $(m,n,M,N)$ with $(m+s,n+t,M+S,N+T)$) to get the $V(\phi(\bar{I}))_K=Z\cup Y$ as described in the lemma so that $Z\subset\phi'(V(\bar{I})_K)$ is a $D$-semianalytic analytic $d$-manifold and $Y$ is the zeroset of an ideal $I'$ of dimension less than $d$.

Next we observe that Theorem 2 of (II) \S III.11.2 of \cite{serre} implies that the set $\pi({\phi'}^{-1}(Z)\cap U)\subset X$ is also a $D$-semianalytic analytic $d$-manifold and that g-dim ${\phi'}^{-1}(Y)<d$ again by Theorem \ref{alleq0} as it is the zero set of the ideal $\phi^{-1}(I')$ which is of dimension less than $d$. Hence we have g-dim $\pi({\phi'}^{-1}(Y))<d$, and therefore we only need to show that $\pi({\phi'}^{-1}(Z))$ satisfies the condition about fibers of points in the parameter space before we are done by induction. 

Let $\pb\in(K^\circ)^{m}\times(K^{\circ\circ})^{n}$, be such that the specialization $X(\pb)\neq\emptyset$, then there is a unique $\qb\in (K^\circ)^{m+s}\times(K^{\circ\circ})^{n+t}$ such that $\pi|_{(K^\circ)^{m+s}\times(K^{\circ\circ})^{n+t}}(\qb)=\pb$. Because $\phi$ respects parameters, $\phi'(\qb)$ is also in the space $(K^\circ)^{m+s}\times(K^{\circ\circ})^{n+t}$ and by Lemma \ref{parajac}, the specialization $(Z\cap U)(\phi'(\qb))$ is either empty or a $D$-semianalytic manifold. Again by the observation that $\phi'$ takes parameters to parameters and fibers to fibers and by Theorem 2 of (II) \S III.11.2 of \cite{serre} we have the statement of the theorem.
\end{proof}

Next we turn our attention to subanalytic sets. We would like to remind the reader that these are projections of D-semianalytic sets onto coordinate hyper-planes and in the case that $K$ is algebraically closed, they are the same as the $D$-semianalytic sets by the Quantifier Elimination Theorem of \cite{rigid}. As a corollary to Theorem \ref{alleq0} we see that when {\em Char $K=0$}, if a $D$-semianalytic set is somewhere dense then it contains an interior point. We are going to prove that this is also true for subanalytic sets by an argument similar to the one we used to prove the Parameterized Smooth stratification Theorem. But first we are going to justify the need for an extra proof by showing that the class of subanalytic sets is a strictly larger class of sets than $D$-semianalytic sets in general through the next easy example.

\begin{ex}\label{subnonsemi}
Let $K=\mathbb{Q}((t))$ with the $t$-adic valuation. Then the set 
$$P_2:=\{x\in K^\circ:(\exists y)(y^2=x)\}$$
is subanalytic but not $D$-semianalytic as it can not be written as a finite boolean combination of discs (see Definition \ref{disc}), contradicting Theorem \ref{complexity}, below.
\end{ex}

Although the previous example shows that in general, unlike a $D$-semianalytic set, {\em complexity} (see Definition \ref{comp}) is not well defined for subanalytic sets, the next theorem shows that subanalytic sets still share a nice geometric property with $D$-semianalytic sets.

\begin{thm}\label{densesub}
If {\em Char $K=0$} and $X$ is a somewhere dense subanalytic subset of $(K^\circ)^m\times(K^{\circ\circ})^n$, then $X$ has an interior point.
\end{thm} 
\begin{proof}
Let  $X'$ be a $D$-semianalytic subset of $(K^\circ)^{m+M}\times(K^{\circ\circ})^{n+N}$ whose projection onto $(K^\circ)^m\times(K^{\circ\circ})^n$ is $X$. Once again we will treat $(K^\circ)^m\times(K^{\circ\circ})^n$ as the parameter space and will be interested in the semianalytic set $X''\subset(K^\circ)^{m+s+M+S}\times(K^{\circ\circ})^{n++tN+T}$ whose image is $X'$ under the coordinate projection 
$$\pi:(K^\circ)^{m+s+M+S}\times(K^{\circ\circ})^{n+t+N+T}\rightarrow(K^\circ)^{m+M}\times(K^{\circ\circ})^{n+N}.$$
Two other coordinate projection maps
$$\pi_1:(K^\circ)^{m+s+M+S}\times(K^{\circ\circ})^{n+t+N+T}\rightarrow(K^\circ)^{m+s}\times(K^{\circ\circ})^{n+t},$$
and 
$$\pi_2:(K^\circ)^{m+s}\times(K^{\circ\circ})^{n+t}\rightarrow (K^\circ)^m\times(K^{\circ\circ})^n$$
will also help us understand the underlying geometry in our setting. The main idea of the proof is to show that $\pi_1(X'')$ contains an open subset of an analytic $(m+n)$-manifold whose charts are given by the coordinate projection $\pi_2$.

 Note that if a finite union of sets is somewhere dense then one of the individual sets must be somewhere dense. Therefore we may assume that  $X'=\Dom B\cap V(I)_K$ for some generalized ring of fractions $B$ over $\Smn$ and ideal $I$ satisfying $I=\ideal(\Dom B\cap V(I)_K)$. By the same fact and Lemma \ref{morenorm} we may also assume that 
$$S_{m'+M',n'+N'}\subset B/\phi(I)$$
is a finite inclusion where $\phi$ is a Weierstrass change of variables respecting parameters. By Lemma 4.4 of \cite{dimtheory} we may assume that $I$ is prime and by Remark \ref{respardense}(ii) we may also assume that $m'+n'=m+n$. 

Write $B=S_{m+s+M+S,n+t+N+T}/J$ for some ideal 
$$J=(\{g_ix_{m+i}-f_i\}_{i=1}^s\cup\{g'_j\rho_{n+j}-f'_j\}_{j=1}^t\cup\{G_ky_{M+k}-F_k\}_{k=1}^S\cup\{G'_l\lambda'_{N+l}-F'_l\}_{l=1}^T),$$ 
and let $\bar{I}\supset J$ be the ideal corresponding to $I$ in $S_{m+s+M+S,n+t+N+T}$ so that $X=\pi_2\circ\pi_1(V(\bar{I})_K\cap U)$, where $U$ is the open set 
\begin{multline*}
U:=\{\pb\in (K^\circ)^{m+s+M+S}\times(K^{\circ\circ})^{n+t+N+T}:\\
\prod_{i=1}^sg_i(\pb)\cdot\prod_{j=1}^tg'_j(\pb)\cdot\prod_{k=1}^SG_k(\pb)\cdot\prod_{l=1}^TG'_l(\pb)\neq 0\}.
\end{multline*}

Observe that $\pi_1(V(I')_K\cap U)$ is contained in the $m+n$ dimensional manifold 
\begin{multline*}
R:=\{\pb\in(K^\circ)^{m+s}\times(K^{\circ\circ})^{n+t}:\pb\in V(\{g_ix_{m+i}-f_i\}_{i=1}^s\cup\{g'_j\rho_{n+j}-f'_j\}_{j=1}^t)\text{ and }\\\prod_{i=1}^sg_i(\pb)\cdot\prod_{j=1}^tg'_j(\pb)\neq 0\},
\end{multline*}
whose charts are given by the restrictions of the projection $\pi_2$ to open subsets of $R$. That is for each $\pb\in R$, there is a rational open neighborhood $W$ of $\pb$, such that $W\cap R$ is given by the relations of the form
\begin{equation}\begin{array}{rcl}\label{chart1}\tag{$\dag$}
x_{m+i}&=&\alpha_{m+i}(x',\rho')\\
\rho_{n+j}&=&\beta_{n+j}(x',\rho')
\end{array}
\end{equation}
where $x'=(x_1,...,x_m)$, $\rho'=(\rho_1,...,\rho_n)$ and $\alpha_{m+i}, \beta_{n+j}\in\mathcal{O}(\pi_2(W))_K$ .

Let us write $\bar{x}, \bar{\rho}, \bar{y}, \bar{\lambda}$ for the images of the variables $x, \rho, y, \lambda$ under the map $\phi$ and observe that k-dim $S_{m+s+M+S,n+t+N+T}/\phi(\bar{I})=m'+n'+M'+N'$. Applying Lemma \ref{parajac} to $V(\phi(\bar{I}))$, we get $Z, Y$ such that $V(\phi(\bar{I}))=Z\cup Y$, where $Z$ is the $(m'+n'+M'+N')$-dimensional manifold which is locally given by the relations of the form 
$$\begin{array}{rcl}
\bar{x}_{m'+i}&=&\bar{\alpha}_{m'+i}(\bar{x}',\bar{\rho}')\\
\bar{\rho}_{n'+j}&=&\bar{\beta}_{n'+j}(\bar{x}',\bar{\rho}')\\
\bar{y}_{M'+k}&=&\bar{\gamma}_{M'+k}(\bar{x}',\bar{\rho}',\bar{y}',\bar{\lambda}')\\
\bar{\lambda}_{N'+l}&=&\bar{\delta}_{N'+l}(\bar{x}',\bar{\rho}',\bar{y}',\bar{\lambda}'),
\end{array}$$ 
where $\bar{x}'=(\bar{x}_1,...,\bar{x}_{m'})$, $\bar{\rho}'=(\bar{\rho}_1,...,\bar{\rho}_{n'})$, $\bar{y}'=(\bar{y}_1,...,\bar{y}_{M'})$, $\bar{\lambda}'=(\bar{\lambda}_1,...,\bar{\lambda}_{N'})$, and g-dim $Y<m'+n'+M'+N'$. Let 
$$\phi':(K^\circ)^{m+s+M+S}\times(K^{\circ\circ})^{n+t+N+T}\rightarrow (K^\circ)^{m+s+M+S}\times(K^{\circ\circ})^{n+t+N+T}$$ 
be the one-to-one analytic map that is induced by $\phi$ so that $\phi'(V(\bar{I})_K)=V(\phi(\bar{I}))_K$, and let $\Delta\not\in\phi(\bar{I})$ be the determinant of the Jacobian as in the proof of Lemma \ref{parajac} so that $Z=\{\pb\in V(\phi(\bar{I})):\Delta(\pb)\neq 0\}$, then notice that $Z\cap \phi'(U)\neq\emptyset$ as otherwise we have 
$$\phi^{-1}(\Delta) \cdot\prod_{i=1}^sg_i\cdot \prod_{j=1}^tg'_j\cdot\prod_{k=1}^SG_k\cdot \prod_{l=1}^TG'_l\in \bar{I}$$ 
contradicting the fact that $\bar{I}$ is prime and $\bar{I}=\ideal(V(\bar{I})_K)$. 

Observe that by Lemma \ref{parajac} we have that $\pi_1(Z)$ (and also $\pi_1(Z\cap\phi'(U))$) is an $(m+n)$-dimensional manifold of $(K^\circ)^{m+s}\times(K^{\circ\circ})^{n+t}$, locally given by the relations
\begin{equation}\begin{array}{rcl}\label{chart2}\tag{$\ddag$}
\bar{x}_{m'+i}&=&\bar{\alpha}_{m'+i}(\bar{x}',\bar{\rho}')\\
\bar{\rho}_{n'+j}&=&\bar{\beta}_{n'+j}(\bar{x}',\bar{\rho}').
\end{array}
\end{equation}
 
Let us write $\psi$ for the Weierstrass change of variables $\phi|_{S_{m+s,n+t}}$ and let $\psi':(K^\circ)^{m+s}\times(K^{\circ\circ})^{n+t}\rightarrow (K^\circ)^{m+s}\times(K^{\circ\circ})^{n+t}$ be the non-singular analytic map that it induces which is also the same as $\phi'|_{(K^\circ)^{m+s}\times(K^{\circ\circ})^{n+t}}$. Let $W_1$ be a rational subdomain of ${(K^\circ)^{m+s}\times(K^{\circ\circ})^{n+t}}$ such that $W_1\cap \pi_1(Z\cap\phi'(U))$ is given by relations of the form \eqref{chart2}. Let 
$$I_1:=\ideal(W_1\cap \pi_1(Z\cap\phi'(U)))\subset\mathcal{O}(W_1)_K$$
so that k-dim $\mathcal{O}(W_1)/I_1=m+n$. Similarly let $W_2\subset{\psi'}^{-1}(W_1)$ be another rational subdomain such that $W_2\cap R$ is given by the relations of the form \eqref{chart1} and let 
$$I_2:=\ideal(R\cap W_2)\subset\mathcal{O}(W_2)_K$$ 
so that k-dim $\mathcal{O}(W_2)/I_2=m+n$ and $\{x_{m+i}-\alpha_{m+i}\}_{i=1}^s\cup \{\rho_{n+j}-\beta_{n+j}\}_{j=1}^t\subset I_2$.  Notice that each $x_{m+i}-\alpha_{m+i}(x',\rho')$ and $\rho_{n+j}-\beta_{n+j}(x',\rho')$ is regular of degree one in $x_{m+i}$ and $\rho_{n+j}$ respectively for each $i$ and $j$. By replacing $W_2$ with a smaller rational open subset and by Lemma \ref{local} we can assume $\mathcal{O}(W_2)_K=T_{m+s+n+t}(\delta,\varepsilon)$ for some non-zero $\delta,\varepsilon\in K^\circ$, and by replacing $W_1$ with a smaller rational open subdomain we may assume ${\psi'}^{-1}(W_1)=W_2$.

 Notice also that in this case $I_2$ is a prime ideal of  $\mathcal{O}(W_2)_K$ by the Weierstrass Division Theorem and observe that ${\psi'}^{-1}$ takes $\pi_1(Z\cap\phi'(U))$ into $R$, hence we have $\psi^{-1}(I_1)\supset I_2$.  Now because $\psi$ induces an isomorphism $\mathcal{O}(W_1)_K\rightarrow \mathcal{O}(W_2)_K$, we have k-dim $\mathcal{O}(W_2)_K/\psi^{-1}(I_1)=m+n$ and therefore $\psi^{-1}(I_1)=I_2$. Hence ${\psi'}^{-1}(V(I_1)_K\cap W_1)$ is the same as the set $V(I_2)_K\cap W_2$ which has an interior point when projected onto $(K^\circ)^{m}\times(K^{\circ\circ})^{n}$ and which is contained in the projection $\pi_1(V(I')_K\cap U)\subset X$.
\end{proof}

As an immediate corollary we have the following result.

\begin{cor}
For a subanalytic subset $X$ of $(K^\circ)^{m}\times(K^{\circ\circ})^{n}$, we have
$$\text{g-dim } X=\text{ w-dim } X.$$
\end{cor}

\section{Piece Number}
In \cite{bart} Bartenwerfer introduced the notions of {\em dimensional filterings} and the {\em piece number of a dimensional filtering} of analytic subvarieties of $(K^\circ)^m\times(K^{\circ\circ})^n$ to generalize the number of Zariski irreducible components of those sets. He showed that for such a set $X\subset (K^\circ)^{m+M}\times(K^{\circ\circ})^{n+N}$ where {\em Char $K=0$} or if {\em Char $K=p>0$} and $[K:K^p]<\infty$, there is a bound $\Gamma$ such that for all $\pb\in (K^\circ)^m\times(K^{\circ\circ})^n$ the fiber $X(\pb)$ has a dimension filtering with piece number less than $\Gamma$. In this section we show that the Parameterized Normalization Lemma (Lemma \ref{parnorm1}) can be used to extend Bartenwerfer's results from analytic varieties to $D$-semianalytic sets in the case {\em Char $K=0$}. 

The concepts of dimensional filterings and piece number can easily be extended to $D$-semianalytic sets, but before we give the  corresponding definitions, we will repeat an example of Bartenwerfer to justify the interest in the piece number rather than the Zariski irreducible components.

\begin{ex}\label{pieceex}
Let $K$ be a discretely valued non-Archimedean complete field of any characteristic with prime element $\varpi$ and let $f=y^2-(\varpi+x)x^2$, then the reduction of $f$ to $(K^\circ/K^{\circ\circ})\left[x,y\right]$ is irreducible and therefore $f$ is itself irreducible. Notice that there is an infinite sequence of points $\pb_i\in V(f)_K$ that converge to the point $(1,1)$ and therefore g-dim $V(f)_K=1$ and $(f)=\ideal(V(f)_K)$. Nevertheless the point $(0,0)$ is an isolated point of $V(I)_K$ as $\varpi$ has no square root in $K$.
\end{ex}

This example shows that the number of irreducible components of an analytic variety may be inadequate as a measure of number of ``pieces'' in that variety. Bartenwerfer introduced the {\em piece number} to investigate those pieces in \cite{bart}. Next we extend the definitions in \cite{bart} to $D$-semianalytic sets.

\begin{defn}\label{dimfilter}
Let $X$ be a $D$-semianalytic subset of $(K^\circ)^m\times(K^{\circ\circ})^n$, with g-dim $X=d\geq 0$, then a $d$-tuple of sets $\mathcal{S}=(S_d,...,S_0)$ is called a dimensional filtering of $X$ if each $S_i$ is a (not necessarily disjoint) finite union of $i$-dimensional $D$-semianalytic $K$-analytic manifolds and $X$ is equal to the  (again not necessarily disjoint) union of the $S_i$.
\end{defn}

In order to be able to define the piece number we also need to define the number of irreducible components of a $D$-semianalytic set.

\begin{defn}\label{irredcomp}
Let $X$ be a $D$-semianalytic subset of $(K^\circ)^m\times(K^{\circ\circ})^n$, let
$$\mathcal{P}:X=\bigcup_{i=1}^r (\dom B_i\cap V(I_i)_K)$$
be a presentation of $X$, let $J_i=\ideal(\dom B_i\cap V(I_i)_K)$, and let $J_i=\cap_{j=1}^{s_i} \mfp_{ij}$ be the irredundant prime decomposition of $J_i$ in $B_i$ for $1\leq i\leq r$, then the number of irreducible components  of this presentation $\mathcal{P}$ of $X$ is defined to be 
$$\#\text{ic }\mathcal{P}:=s_1+...+s_r,$$
and the number of irreducible components of $X$ is defined to be
$$\#\text{ic }X:=\inf \{\#\text{ic }\mathcal{P}:\mathcal{P}\text{ is a presentation of }X\}.$$
\end{defn}

Now we are ready to define

\begin{defn}\label{piecen}
Given a $D$-semianalytic set $X$ with dimensional filtering $\mathcal{S}=(S_d,...,S_0)$, we define the piece number of $\mathcal{S}$ as
$$\text{pn }\mathcal{S}:=\sum_{i=0}^d (\#\text{ic }S_i).$$
We will define the piece number of $X$ to be 
$$\text{pn }X:= \inf\{\text{pn }\mathcal{S}:\mathcal{S}\text{ is a dimensional filtering of }X\}.$$
\end{defn}

Note that this definition is slightly different from Bartenwerfer's definition as instead of working with {\em pure dimensional} subsets of varieties we work with analytic manifolds. However, as in his definition, pn $X$ is easily seen to dominate the number of irreducible components and isolated points in $X$. It is also easy to see that the analytic variety in Example \ref{pieceex} has piece number at least two, confirming what one would intuitively expect from the piece number.

The following theorem is the generalization of the main theorem of \cite{bart} which was proved for analytic varieties in the case {\em Char }$K=0$ or {\em Char $K=p>0$} and $[K:K^p]<\infty$. As we usually do throughout this paper, we will assume that {\em Char $K=0$}.

\begin{thm}\label{boundedpn}
Let {\em Char $K=0$} and let $X\subset(K^\circ)^{m+M}\times(K^{\circ\circ})^{n+N}$ be a $D$-semianalytic set, then there is a bound $\Gamma\in\mathbb{N}$ such that for all $\pb\in(K^\circ)^{m}\times(K^{\circ\circ})^{n}$ the piece number of the fiber $X(\pb)$ is less than $\Gamma$.
\end{thm}
\begin{proof}
By definition it is enough to find a $\Gamma$ such that for each $\pb\in(K^\circ)^{m}\times(K^{\circ\circ})^{n}$, $X(\pb)$ has a dimensional filtering with piece number less than $\Gamma$. Our plan is to make use of the Parameterized Normalization Lemma (Lemma \ref{parnorm1}) to first normalize the $D$-semianalytic set we are working on and then use Parameterized Smooth Stratification Theorem (Theorem \ref{parass}) to find a dimensional filtering for each fiber. The piece number of each such a dimensional filtering will be uniformly bounded by the product of degrees of minimal polynomials of the integral variables of the normalization.

Let $X_1,...,X_k$ be as in Theorem \ref{parass} so that $X=X_1\cup...\cup X_k$ and for all $\pb\in(K^\circ)^{m}\times(K^{\circ\circ})^{n}$ $X_i(\pb)$ is either empty or a $D$-semianalytic manifold. We will make use of the construction in the proof of Theorem \ref{parass} to further observe that in this case for each $i$ we have
$$X_i=\Dom B_i\cap V(\mfp_i)_K\cap\{\qb\in(K^\circ)^{m+M}\times(K^{\circ\circ})^{n+N}:\Delta_i(\qb)\neq 0\}$$
for some generalized ring of fractions $B_i$ over $\Smn$, prime ideal $\mfp_i\subset B_i$ satisfying $\mfp_i=\ideal(\Dom B_i\cap V(\mfp_i)_K)$ and determinants of Jacobians $\Delta_i\in B_i$. Moreover there are integers $m_i,M_i,n_i, N_i$ and Weierstrass changes of variables $\phi_i$ such that there is a finite injection
$$S_{m_i+M_i,n_i+N_i}\rightarrow B_i/\phi_i(\mfp_i).$$

Note that if $z$ is a variable not appearing in the presentation of any $B_i$ then $X_i$ is a $D$-semianalytic set with the presentation
$$X_i=\Dom (B_i\ls z\rs/\Delta_i z-0)\cap V(\mfp_i\cdot (B_i\ls z\rs/\Delta_iz-0))_K.$$

Now for a given $\pb\in(K^\circ)^{m}\times(K^{\circ\circ})^{n}$ let g-dim $X(\pb)$ be $d$ and set
$$S_j(\pb):=\bigcup\{X_i: X_i(\pb)\text{ is }j\text{-dimensional}\}$$
for all $0\leq j\leq d$ so that $\mathcal{S}:=(S_d(\pb),...,S_0(\pb))$ is a dimensional filtering for $X(\pb)$. Now it is enough to show that for each $X_i$ there is a bound $\Gamma$ such that for $\pb\in(K^\circ)^{m}\times(K^{\circ\circ})^{n}$ we have $\#$ic $X_i(\pb)\leq \Gamma$. 

Let $B_i=S_{m+s_i+M+S_i,n+t_i+N+T_i}/I_i$ and let $\bar{\mfp}_i$ be the ideal corresponding to $\mfp$ in the quasi-affinoid algebra $S_{m+s_i+M+S_i,n+t_i+N+T_i}/I_i$. By $\psi_i$, let us denote the restriction of the Weierstrass change of variables $\phi_i$ to $S_{m+s_i,n+t_i}$ so that there is a finite homomorphism
$$S_{m_i,n_i}\rightarrow S_{m+s_i,n+t_i}/\psi_i(\bar{\mfp}_i\cap S_{m+s_i,n+t_i}).$$ 
For $\pb\in(K^\circ)^{m}\times(K^{\circ\circ})^{n}$ let $B_i(\pb)$ and $\mfp_i(\pb)$, $\Delta_i(\pb)$ denote the objects we obtain by making the obvious substitutions the notation indicates so that $B_i(\pb)$ is a generalized ring of fractions over $S_{M,N}$, $\Delta_i(\pb)\in B_i(\pb)$ and $\mfp_i(\pb)$ is an ideal of $B_i(\pb)$.

Note that the restriction $\phi_i^*$ of $\phi_i$ to $y$ and $\lambda$ variables induces a finite homomorphism
$$S_{M_i,N_i}\rightarrow B_i(\pb)/\phi_i^*(\mfp_i(\pb)),$$
and if $X_i(\pb)$ is not empty, then $X_i(\pb)$ is an $(M_i+N_i)$-dimensional manifold and it is easy to see also that g-dim $X_i(\pb)=M_i+N_i$ and therefore k-dim $B_i(\pb)/\phi_i^*(\mfp_i(\pb))\geq M_i+N_i$ and we see that in fact the above homomorphism is an injection. 

Now put 
$$J_{\pb}:=\ideal(\text{Dom}_{M+N}(B_i(\pb)\ls z\rs/\Delta_i(\pb)z-0)\cap V(\mfp_i(\pb)\cdot(B_i(\pb)\ls z\rs/\Delta_i(\pb)z-0))_K)$$ 
and let us write $\bar{J}_{\pb}$ for the corresponding ideal in $S_{M+S_i+1,N+T_i}$. Notice that $z\in\bar{J}_{\pb}$ by the definition of $J_{\pb}$. Therefore $(B_i(\pb)\ls z\rs/\Delta_i(\pb)z-0)/\phi_i^*(J_{\pb})$ is finite over $S_{M_i,N_i}$. Notice also that $X_i(\pb)=\text{Dom}_{M,N}(B_i(\pb)\ls z\rs/\Delta_i(\pb)z-0)\cap V(J_{\pb})_K$ and therefore g-dim $X_i(\pb)=M_i+N_i$. In fact we have a finite inclusion
$$S_{M_i,N_i}\subset (B_i(\pb)\ls z\rs/\Delta_i(\pb)z-0)/\phi_i^*(J_{\pb}).$$

Note that over each prime ideal of $S_{M_i,N_i}$ there are at most 
$$\Gamma_{\pb}=r_1\cdots r_{(M+S_i-M_i)+(N+T_i-N_i)}$$ 
--the product of degrees of the lowest degree monic polynomials in $\phi_i^*(\bar{J}_{\pb})$ of the variables $y_{M_i+1},...,y_{M+S_i}$ and $\lambda_{N_i+1},...,\lambda_{N+T_i}$ over $S_{M_i,N_i}$-- many prime ideals of $(B_i(\pb)\ls z\rs/\Delta_i(\pb)z-0)/J_{\pb}$. Now the statement follows from the fact that for each $\pb\in(K^\circ)^m\times(K^{\circ\circ})^n$,  $\Gamma_{\pb}$ is dominated by 
$$\Gamma=s_1\cdots s_{(m+s_i-m_i)+(M+S_i-M_i)+(n+t_i-n_i)+(N+T_i-N_i)}$$ 
--the product of degrees of the lowest degree monic polynomials in $\phi_i(\bar{\mfp}_i)$ of the variables $x_{m_i+1},...,x_{m+s_i}$, $\rho_{n_i+1},...,\rho_{n+t_i}$, $y_{M_i+1},...,y_{M+S_i}$ and $\lambda_{N_i+1},...,\lambda_{N+T_i}$ over $S_{m_i+M_i,n_i+N_i}$.
\end{proof}

\section{Complexity}
Next we will turn our attention to the case that the fibers of a $D$-semianalytic set are subsets of $K^\circ$. Our aim is to give a simpler proof of the main theorem of \cite{complex} using the results of the previous section on the piece number. That is, we are going to prove that given a $D$-semianalytic set $X$, there is a uniform bound $\Gamma$ such that each one dimensional fiber of $X$ is a boolean combination of at most $\Gamma$ discs (see Definition \ref{disc} below) and points, by showing that the piece number of each such fiber is closely related to the number of discs and points required for such a boolean combination. This result is analogous to the main theorem of \cite{onedim} in the sense that both concern a measure of complexity of one-dimensional fibers of definable sets. However, at the end of this section we will prove a more readily recognizable analog of that theorem in our setting.

Before we state and prove this theorem we need some groundwork.

\begin{defn}\label{disc}
A disc in $K^\circ$ is a set of the form
$$\begin{array}{rcl}
D^-(a,r)&:=&\{x\in K^\circ:|x-a|<r\},\text{ or }\\
D^+(a,r)&:=&\{x\in K^\circ:|x-a|\leq r\}
\end{array}$$
where $a\in K^\circ$ and $r\in\sqrt{|K^\circ|}$.

We will follow the terminology of \cite{complex} and call a set of the form a disc minus a finite union of of discs {\em a $K$-rational special set}, or in case $K$ is algebraically closed, a special set. 
\end{defn}

In \cite{linepln} Lipshitz and Robinson showed that for an algebraically closed complete non-archimedean field $\bar{K}$, if $X$ is an $R$-domain in $\bar{K}^\circ$ then $X$ is a finite union of special sets. On the other hand with a little effort one sees that any such union is an $R$-domain. Let $K\subset \bar{K}$ be another non-Archimedean complete field which is not necessarily algebraically closed. As any member of the disc is the center of the disc, it is easy to see that given a disc $C\subset\bar{K}^\circ$, either $C\cap K^\circ=\emptyset$ or $C\cap K^\circ$ is a disc in $K^\circ$. Therefore any $K$-rational $R$-domain of $K$ for an arbitrary non-Archimedean complete field $K$ is a finite union of $K$-rational special sets. Notice that there may be more than one way of writing $X$ as a union of disjoint special sets as in the case the residue field $K^\circ/K^{\circ\circ}$ is finite, the closed unit disc $D^+(0,1)$ is a disjoint union of the open unit disc $D^-(0,1)$ and finitely many smaller closed discs. 

\begin{defn}\label{comp}
Let $X\subset K$ be a $K$-rational $R$-domain and 
$$\mathcal{D}:X=\bigcup_{i=1}^r S_i=\bigcup_{i=1}^r(C_i\setminus\bigcup_{j=1}^{s_i}C_{ij})$$
be a decomposition of $X$ into $K$-rational special sets $S_i=C_i\setminus\bigcup_{j=1}^{s_i}C_{ij}$ then the complexity of $\mathcal{D}$ is defined to be
$$\text{comp }\mathcal{D}:=r+s_1+...+s_r,$$
and the complexity of the $R$-domain $X$ is defined to be
\begin{multline*}
\text{comp } X:=\inf \{\text{comp }\mathcal{D}:\mathcal{D}\text{ is a decomposition of $X$ into}\\
\text{$K$-rational special sets}\}.
\end{multline*}
\end{defn}

Let $X=\dom B$ be a $K$-rational $R$-domain for some generalized ring of fractions $B$ over $S_{0,1}(E,K)$ (or $S_{0,1}(E,K)$) and let $\bar{X}$ be the $R$-domain $\bkdom B$, then it is clear from the definition that the complexity of $\bar{X}$ dominates the complexity of $X$. 

Our aim is to establish a relation between the piece number and the complexity. More specifically we are going to show that  for a $K$-rational $R$-domain $X\subset K^\circ$, 
$$\text{pn $X+$ pn $(K^\circ\setminus X)\geq$ comp }X.$$
For this we need a deeper understanding of some special type of quasi-affinoid algebras. The following two key lemmas which we state without proof are going to be proved in \cite{cell} by Cluckers, Lipshitz and Robinson. Both of these lemmas are well known in the affinoid category and we refer the reader to $\S$7.3.3 of \cite{bgr} and $\S$2.2 of \cite{vdPF} for further details.

\begin{lem}\label{directsum}
Suppose $K$ is an algebraically closed complete field and let $B$ be the ring of analytic functions of an $R$-domain $X\subset K^\circ$. Let 
$$X=\bigcup_{i=1}^r X_i$$
be a decomposition of $X$ into disjoint special sets with minimal complexity, then we have
$$B=\mathcal{O}(X)=\bigoplus_{i=1}^r\mathcal{O}(X_i).$$
\end{lem}

We would like to point out that one can find an analogue of the above lemma in the affinoid category in Proposition 7 of $\S$ 7.3.2 of \cite{bgr}.

\begin{lem}=\label{pid1}
Suppose $K$ is an algebraically closed complete field and suppose $B$ is the ring of analytic functions of a special set, then for any $h\in B$ we have
$$h=u\cdot p(x)$$
where $u\in B$ is a unit and $p(x)\in K\left[ x\right]$.
\end{lem}

As an immediate corollary to the above lemmas we have:

\begin{lem}\label{pid2}
Suppose $K$ is an algebraically closed complete field and suppose $B$ is the ring of analytic functions of an $R$-subdomain of $K^\circ$, then $B$ is a principal ideal domain.
\end{lem}
\begin{proof}
By separating Dom$_{1,0}B$ into disjoint special sets and making use of Lemma \ref{directsum} we may assume that Dom$_{1,0}B$ is a special set. Now it is enough to show that for any $f_1,f_2\in B$ there is an $f\in B$ such that $(f_1,f_2)=(f)$. Using Lemma \ref{pid1} write $f_i=u_i\cdot p_i(x)$  where $u_i\in B$ is a unit and $p_i(x)\in K\left[ x\right]$ for $i=1,2$. Now put $f$ to be the greatest common divisor of $p_1(x)$ and $p_2(x)$ and the result follows.
\end{proof}

The next lemma shows that the an arbitrary generalized ring of fractions over $S_{1,0}$ (or $S_{0,1}$) is not too far removed from one that is the ring of analytic functions of an $R$-domain. But before we state and prove the lemma, we need a notation which enables us to keep track of the inductive steps in the inductive construction of a generalized ring of fractions.

Let $B$ be a generalized ring of fractions over $S_{1,0}$ which is constructed in $m+n$ steps. That means that there is a sequence $\{B_i\}_{i=0}^{m+n}$ of generalized rings of fractions such that $B_0=S_{1,0}$, $B_{m+n}=B$ and if $B_{i}$ is given by the presentation
$$B_{i}=S_{1+m_i,n_i}/(\{g'_jx_j-f'_j\}_{j=1}^{m_i}\cup\{g''_l\rho_l-f''_l\}_{l=1}^{n_i})$$
where $m_i+n_i=i$ then $B_{i+1}$ is given either by $B_{i}\ls x_{m_i+1}\rs/(gx_{m_i+1}-f)$ or by $B_{i}\lb \rho_{n_i+1}\rb/(g\rho_{n_i+1}-f)$ for some $f,g\in S_{1+m_i,n_i}$ for $i=0,...,m+n-1$. For such an inductive construction of $B$, we will introduce the following notation inductively: at the $(i+1)^{st}$ step define
$$\begin{array}{rcl}
\xr^B_{i+1}&:=&\begin{cases}
x_{m_i+1}\text{ if }B_{i+1}=B_{i}\ls x_{m_i+1}\rs/(fx_{m_i+1}-g)\\
\rho_{n_i+1}\text{ if }B_{i+1}=B_{i}\lb \rho_{n_i+1}\rb/(f\rho_{n_i+1}-g)\\
\end{cases}\\
\left[f\right]^B_{i+1}&:=&f\\
\left[g\right]^B_{i+1}&:=&g
\end{array}$$
So that
$$B=S_{1+m,n}/(\{\left[g\right]^B_i\xr^B_i-\left[f\right]^B_i\}_{i=1}^{m+n}).$$
When the generalized ring of fractions $B$ is clear from the context, we will just write $\xr_i$, $f_i$, $g_i$ instead of $\xr^B_i$, $\left[f\right]^B_i$, $\left[g\right]^B_i$. Now we are ready to state and prove:

\begin{lem}\label{Rific}
Suppose $K$ is algebraically closed, complete. Let $B$ be the ring of analytic functions of an $R$-domain $X\subset K^\circ$, $X\neq\emptyset$. Write
$$B=S_{1+m,n}/(\{g_i\xr_i-f_i\}_{i=1}^{m+n}),$$
where $g_i,f_i\in S_{1+m_i,n_i}$, then for each $i=1,...,m+n$ there exist $\bag_i,\baf_i\in S_{1+m_i,n_i}$ for such that

i) $(\bag_1\xr_1-\baf_1,...,\bag_{i-1}\xr_{i-1}-\baf_{i-1},\bag_i,\baf_i)$ is the unit ideal of $S_{1+m_i,n_i}$ and 
$$\sqrt{(\{g_i\xr_i-f_i\}_{i=1}^{m+n})}\subset\sqrt{(\{\bag_i\xr_i-\baf_i\}_{i=1}^{m+n})},$$

ii) for $\bar{B}:=S_{1+m,n}/(\{\bag_i\xr_i-\baf_i\}_{i=1}^{m+n})$, Dom$_{1,0}\bar{B}\setminus X$ consists of finitely many points,

iii) any minimal prime divisor $\mfp\subset S_{1+m,n}$ of $(\{\bag_i\xr_i-\baf_i\}_{i=1}^{m+n})$ is also a minimal prime divisor of $(\{g_i\xr_i-f_i\}_{i=1}^{m+n})$.
\end{lem}
\begin{proof}
We will start by observing that by Lemma \ref{pid1} if $A$ is a generalized ring of fractions over $S_{1,0}$ and Dom$_{1,0} A$ is a special set, then $A$ is an integral domain and k-dim $A=1$. Therefore if we found $\bag_1,\baf_1,...,\bag_{m+n},\baf_{m+n}$ satisfying the conditions (i) and (ii) above, then by Lemma \ref{directsum} k-dim $\bar{B}=1$. Hence for any minimal prime divisor $\mfp$ of $(\{\bag_i\xr_i-\baf_i\}_{i=1}^{m+n})$, k-dim $S_{1+m,n}/\mfp=1$

Next we are going to proceed by induction on $m+n$. The reader may also think of this as an induction on complexity of $\mathcal{L}_{\text{an}}^D$-terms that appear in the definition of Dom$_{m+1,n}B$, where $\mathcal{L}_{\text{an}}^D$ stands for the analytic language of \cite{rigid}. Suppose the lemma holds for all generalized rings of fractions constructed in less than $m+n$ steps. Put
$$B'=S_{1+m_{m+n-1},n_{m+n-1}}/(\{g_i\xr_i-f_i\}_{i=1}^{m+n-1}.$$
That is, let $B'$ be the the ring we obtain in the penultimate step in the construction of $B$, and for simplicity in notation, let us assume that $\xr_{m+n}$ is $\rho_n$ so that $B'=S_{1+m,n-1}/(\{g_i\xr_i-f_i\}_{i=1}^{m+n-1})$ and $B=B'\lb\rho_n\rb/\bar{g}_{m+n}\rho_n-\bar{f}_{m+n}$, where $\bar{g}_{m+n}$ and $\bar{f}_{m+n}$ are the images of $g_{m+n}$ and $f_{m+n}$ in $B'$.

Assume that we found $\bag_1,\baf_1,...,\bag_{m+n-1},\baf_{m+n-1}$ satisfying the conclusion of the lemma for $B'$, then $B''= S_{1+m,n-1}/(\{\bag_i\xr_i-\baf_i\}_{i=1}^{m+n-1})$ is a ring of functions of an $R$-domain and if $U_1,...,U_k$ are the disjoint special sets that make up Dom$_{1,0}B'$ then $B''=\bigoplus_{j=1}^k\mathcal{O}(U_j)$. Notice that k-dim $B''=1$ and therefore  the possible values for
$$d:=\text{k-dim }S_{1+m,n-1}/(\{\bag_i\xr_i-\baf_i\}_{i=1}^{m+n-1}, g_{m+n},f_{m+n})$$
are $-1$, $0$ and $1$.

{\em Case 1:} $d=-1$.

Then we put $\bag_{m+n}=g_{m+n}$, $\baf_{m+n}=f_{m+n}$ and the statements (i) and (ii) are easily seen to be satisfied. Now assume $\mfp\in S_{1+m,n}$ is a minimal prime divisor of $(\{\bag_i\xr_i-\baf_i\}_{i=1}^{m+n})$, by observing the Krull dimensions we see that $\mfp\cap S_{1+m,n-1}$ is minimal over $(\{\bag_i\xr_i-\baf_i\}_{i=1}^{m+n-1}$ and by inductive hypothesis also over $(\{g_i\xr_i-f_i\}_{i=1}^{m+n-1}$. Now if $\mfq\subset\mfp$ is a minimal prime divisor of $(\{g_i\xr_i-f_i\}_{i=1}^{m+n})$ then $\mfq\cap S_{1+m,n-1}$ contains $(\{\bag_i\xr_i-\baf_i\}_{i=1}^{m+n-1})$ and therefore $\mfp=\mfq$.

{\em Case 2:} $d=0$.

By Lemma \ref{pid2} there is an $h\in S_{1+m,n-1}$ such that 
$$(\{\bag_i\xr_i-\baf_i\}_{i=1}^{m+n-1}, g_{m+n},f_{m+n})=(\{\bag_i\xr_i-\baf_i\}_{i=1}^{m+n-1}, h)$$
and the image $\bar{h}$ of $h$ in $B''$ is of the form $\bar{h}=p_1(x)\oplus...\oplus p_k(x)$ where $p_j(x)\in K\left[x\right]$ for all $j$. Write
$$\begin{array}{rcl}
f_{m+n}&=&F_{m+n}h+a\\
g_{m+n}&=&G_{m+n}h+b
\end{array}$$
where $a,b\in (\{\bag_i\xr_i-\baf_i\}_{i=1}^{m+n-1})$. By Nullstellensatz and Lemma \ref{pid1} it is easily seen that $(\{\bag_i\xr_i-\baf_i\}_{i=1}^{m+n-1}, \bag_{m+n},\baf_{m+n})$ is the unit ideal. It is also easy to verify that the conclusion of the lemma holds if we replace $(\{g_i\xr_i-f_i\}_{i=1}^{m+n})$ with $(\{\bag_i\xr_i-\baf_i\}_{i=1}^{m+n-1},g_{m+n}\rho_n-f_{m+n})$ and we are done by {\em Case 1}.

{\em Case 3:} $d=1$.

Let us write $\bar{f}_{m+n}=a_1\oplus...\oplus a_k$ and $\bar{g}_{m+n}=b_1\oplus...\oplus b_k$ for the images of $f_{m+n}$ and $g_{m+n}$ in $B''$. Then after a rearrangement of the components there is an $l\leq k$ such that $a_j=b_j=0_j$ for $1\leq j\leq l$ and k-dim $\mathcal{O}(U_j)/(a_j,b_j)$ is $0$ or $-1$ for $l<j\leq k$. Let $h\in S_{1+m,n-1}$ be such that the image of $h$ in $B''$ is $1_1\oplus...\oplus1_l\oplus0_{l+1}\oplus...\oplus 0_k$ and set $\baf_{m+n}:=f_{m+n}+h$, $\bag_{m+n}:=g_{m+n}$ so that k-dim $S_{1+m,n}/(\{\bag_i\xr_i-\baf_i\}_{i=1}^{m+n})$ is $0$ or $-1$. By Nullstellensatz it is again easy to see that the conclusion of the lemma holds if we replace $(\{g_i\xr_i-f_i\}_{i=1}^{m+n})$ with $(\{\bag_i\xr_i-\baf_i\}_{i=1}^{m+n-1},g_{m+n}\rho_n-f_{m+n})$ and the problem reduces to {\em Case 1} or {\em 2}.
\end{proof}

Now we are ready to give another proof of the main theorem (Theorem 1.6) of \cite{complex}. Note that this result is analogous to the main theorem of \cite{onedim} in that it puts a bound on how ``complicated'' one dimensional fibers of $D$-semianalytic sets can get.

\begin{thm}\label{complexity}
Let {\em Char $K=0$} and let $X\subset (K^\circ)^{m+1}\times (K^{\circ\circ})^n$ be $D$-semianalytic, then there exists a bound $\Gamma$ such that for each $\pb\in(K^\circ)^{m}\times (K^{\circ\circ})^n$ there is a $K$-rational $R$-domain $Y(\pb)\subset K^\circ$ of complexity at most $\Gamma$  and a set $Q(\pb)\subset K^\circ$ of at most $\Gamma$ points such that
$$Y(\pb)\setminus Q(\pb)=X(\pb)\setminus Q(\pb).$$
\end{thm}
\begin{proof}
By the explanation following Definition \ref{disc} we can assume that $K$ is algebraically closed. Let $\Gamma_1\in\mathbb{N}$ be such that pn $X(\pb)<\Gamma_1$ for all $\pb\in (K^\circ)^{m}\times (K^{\circ\circ})^n$ as in Theorem \ref{boundedpn}. Fix a $\pb\in (K^\circ)^{m}\times (K^{\circ\circ})^n$ and let $\mathcal{S}=(S_1,S_0)$ be a dimensional filtering of $X(\pb)$ such that pn $\mathcal{S}=$ pn $X(\pb)$. Write $S_1=\bigcup_{i=1}^k$ Dom$_{1,0}B_i\cap V(I_i)$ for some generalized rings of fractions $B_i$ over $S_{1,0}$ and ideals $I_i\subset B_i$ where each Dom$_{1,0}B_i/I_i$ is a one dimensional $K$-analytic $1$-manifold. 

Let $\bar{B}_i$ denote the generalized ring of fractions we obtain from $B_i$ by applying Lemma \ref{Rific} so that each Dom$_{1,0}\bar{B}_i$ is an $R$-domain and by Lemma \ref{directsum} so is each Dom$_{1,0}\bar{B}_i\cap V(I_i)$ with $\#$ic Dom$_{1,0}\bar{B}_i\cap V(I_i)$--many special sets in any decomposition. Define $Y(\pb):=\bigcup_{i=1}^k$ Dom$_{1,0}\bar{B}_i\cap V(I_i)$ and $Q(\pb):=(Y(\pb)\setminus X(\pb))\cup (X(\pb)\setminus Y(\pb))$, so that $Y(\pb)$ is an $R$-domain and number of special sets in any decomposition of $Y(\pb)$ is bounded by
$$s:=\sum_{i=1}^k\#\text{ic }(\text{Dom}_{1,0} B_i\cap V(I_i))\geq\sum_{i=1}^k\#\text{ic }(\text{Dom}_{1,0} \bar{B}_i\cap V(I_i)).$$

Notice that $s+|X(\pb)\setminus Y(\pb)|\leq s+ |Q(\pb)|\leq\Gamma_1$, and $S_0\cap(Y(\pb)\setminus S_1)=\emptyset$ as if $\qb\in S_0\cap(Y(\pb)\setminus S_1)$ then using the construction in Lemma \ref{Rific}, it is easy to see that $(\{\qb\}\cup S_1,S_o\setminus\{\qb\})$ is another dimensional filtering of $X(\pb)$ with a smaller piece number. Therefore the points in $Y(\pb)\setminus S_1$ are isolated points of $K^\circ\setminus X(\pb)$. Hence the bound $\Gamma_2$ for the piece numbers of fibers  of the $D$-semianalytic set $((K^\circ)^{m+1}\times (K^{\circ\circ})^n)\setminus X$ dominates the total number of holes in the special sets that make up $X(\pb)$ plus the number of points in $Y(\pb)\setminus S_1$ and $\Gamma:=\Gamma_1+\Gamma_2$ is the desired bound in the statement of the lemma.
\end{proof}

As pointed out by Leonard Lipshitz we can prove a theorem with a statement more similar to Theorem A of \cite{onedim} in all characteristics using an argument similar to the one that is used to prove the Parameterized Normalization Lemma. Before we state theorem we would like to remind the reader that a {\em semialgebraic subset} of  $(K^\circ)^m\times(K^{\circ\circ})^n$ is a finite union of sets of the form
$$\{\pb\in(K^\circ)^m\times(K^{\circ\circ})^n:\bigwedge_{i=1}^s(|f_i(\pb)|<|g_i(\pb)|)\wedge\bigwedge_{j=1}^t(|f'_j(\pb)|\leq|g'_j(\pb)|)\}$$
where $f_i$, $g_i$, $f'_j$, and $g'_j$ are polynomials.

\begin{thm}\label{thmA}
Let $X\subset (K^\circ)^{m+1}\times (K^{\circ\circ})^n$ be a $D$-semianalytic set, then there exists a {\em semialgebraic set} $Y\subset (K^\circ)^{m+M+1}\times (K^{\circ\circ})^{n+N}$ such that for every $\pb\in(K^\circ)^{m+1}\times (K^{\circ\circ})^n$ there is a $\qb\in(K^\circ)^{m+M}\times (K^{\circ\circ})^{n+N}$ such that the fibers $X(\pb)$ and $Y(\qb)$ are the same.
\end{thm}
\begin{proof}
For the moment we will assume that $X=\text{Dom}_{m+1,n}B\cap V(I)_K$ for some generalized ring of fractions $B$ over $S_{m+1,n}$ and ideal $I\subset B$. Write 
$$B=S_{m+1+s,n+t}/(\{g_ix_i-f_i\}_{i=m+1}^{m+s}\cup\{g'_j\rho_j-f'_j\}_{j=n+1}^{n+t}),$$
where $S_{m+1+s,n+t}$ denotes the ring of separated power series over the variables $x_1,...,x_m,y,x_{m+1},...,x_{m+s}$, $\rho_1,...,\rho_{n+t}$ and $y$ is the variable corresponding to the fiber spaces. Let $\bar{I}$ be the ideal corresponding to $I$ in $S_{m+1+s,n+t}$ and let $h_1,...,h_k$ generate $\bar{I}$. For $i=1,...,k$ write
$$h_i=\sum_{j\in \mathbb{N}}a_{i,j}(x,\rho)y^j,$$
let $\bar{I}'$ be the ideal of $S_{m+s,n+t}$ generated by $\{a_{i,j}\}_{i,j}$, and let $I'$ be the corresponding ideal in $B$. By Lemma 3.1.6 of \cite{RoSPS} there is a finite set $Z\subset\{1,...,k\}\times\mathbb{N}$ such that for any $\beta\in \{1,...,k\}\times\mathbb{N}$, there are $\{b_{\beta,\alpha}\}_{\alpha\in Z}$ such that $a_\beta=\sum_{\alpha\in Z}b_{\beta\alpha}a_{\alpha}$ and $||a_{\beta}||\geq||b_{\beta,\alpha}a_\alpha||$ for all $\alpha\in Z$.

Now for each $\alpha=(i_0,j_0)\in Z$ define
$$B_\alpha=B\ls \left\{a_{i_0,j}/a_{i_0,j_0}\right\}_{(i_0,j)\in Z,j<j_0}\rs\lb\left\{a_{i_0,j}/a_{i_0,j_0}\right\}_{(i_0,j)\in Z, j>j_0}\rb,$$
$X_\alpha=\text{Dom}_{m+1,n}B_\alpha\cap V(IB_\alpha)_K$ and $X'=\text{Dom}_{m+1,n}B\cap V(I')_K$. Notice that 
$$X=X'\cup\bigcup_{\alpha\in Z}X_\alpha.$$

Next we will show the statement of the theorem holds if we replace $X$ with $X_\alpha$ or $X'$. For $\pb\in(K^\circ)^{m}\times (K^{\circ\circ})^{n}$ the fiber $X'(\pb)$ is either empty or $K^\circ$ and so the statement holds trivially for $X'$. On the other hand if we write $B_\alpha=S_{m+1+s_\alpha,n+t_\alpha}/J_\alpha$ and $\bar{I}_\alpha$ for the ideal of $S_{m+1+s_\alpha,n+t_\alpha}$ corresponding to $IB_\alpha$ then for $\alpha=(i_0,j_0)\in Z$ we have $h_{i_0}\in\bar{I}_{i_0,j_0}$ and we can write
$$h_{i_0}\equiv a_\alpha(y^{j_0}+\sum_{j\neq j_0} b_{i_0,j} y^j)\text{ mod }(\{a_\alpha x_l-a_{i_0,j}\}_{l=s+1}^{s_\alpha}\cup \{a_\alpha\rho_l-a_{i_0,j}\}_{l=t+1}^{t_\alpha})$$
for some $b_{i_0,j}\in S_{m+s_\alpha,n+t_\alpha}$ which make $y^{j_0}+\sum_{j\neq j_0} b_{i_0,j} y^j$ regular in $y$ of degree $j_0$. Therefore by the Weierstrass Preparation theorem we see that there is a monic regular polynomial $p_\alpha(y)$ of degree $j_0$ in $S_{m+s_\alpha,n+t_\alpha}\left[ y\right]$ which, when treated as a function over $(K^\circ)^{m+1}\times(K^{\circ\circ})^n$ vanishes at $X_\alpha$. Let $\bar{p}_\alpha(y)$ be the image of $p_\alpha(y)$ in $B_\alpha$ and write $\bar{p}_\alpha(y)=y^{j_0}+c_{\alpha,j_0-1}y^{j_0-1}+...+c_{\alpha,0}$ for some $D$-functions $c_{\alpha,j_0-1},...,c_{\alpha,0}$ over $(K^\circ)^m\times(K^{\circ\circ})^n$. 

Now let $\bar{g}_i$, $\bar{f_i}$, $\bar{g}'_j$, $\bar{f}'_j$, $\bar{a}_{i_0,j}$, $\bar{h}_i$ denote the images of remainders of $g_i$, $f_i$, $g'_j$, $f'_j$, $a_{i_0,j}$, $h_i$ after the application of Weierstrass Division Theorem to divide by $p_\alpha(y)$ in the generalized ring of fractions $B_\alpha$. Then $X_\alpha$ consists of $\pb\in(K^\circ)^{m+1}\times(K^{\circ\circ})^n$ such that
$$\begin{array}{rcll}
|\bar{f}_i(\pb)|&\leq&|\bar{g}_i(\pb)|\neq 0&\text{ for }i=m+1,...,m+s\\
|\bar{f}'_j(\pb)|&<&|\bar{g}'_j(\pb)|&\text{ for }j=n+1,...,n+t\\
|\bar{a}_{i_0,j}(\pb)|&\leq&|\bar{a}_{i_0,j_0}(\pb)|\neq 0&\text{ for } (i_0,j)\in Z, j<j_0\\
|\bar{a}_{i_0,j}(\pb)|&<&|\bar{a}_{i_0,j_0}(\pb)|&\text{ for } (i_0,j)\in Z, j>j_0\\
\bar{h}_l(\pb)&=&0&\text{ for }l=1,...,k\\
\bar{p}_\alpha(y)&=&0.
\end{array}$$

On each of the lines above $y$ appears polynomially, and therefore by introducing new variables for each of the non-polynomial ($D$-function) terms in the above description we get a semialgebraic set $Y_\alpha$. Given a $\pb\in(K^\circ)^{m+1}\times(K^{\circ\circ})^n$, if we substitute the values coming from those $D$-functions of $\pb$ for those new variables the resulting fiber is the same as $X_\alpha(\pb)$ and we have the statement of the theorem for $X_\alpha$.

On the other hand, by introducing different variables for each $X_\alpha$ in the above process we can get the fiber $X(\pb)$ as the union of fibers $Y_{\alpha}(\pb)$. As a $D$-semianalytic set is a finite union of sets of the form $\text{Dom}_{m+1,n}B\cap V(I)_K$, this argument proves the theorem.
\end{proof}

\begin{rem}
We would like to note that the above proof can easily be modified to actually give another proof of Theorem A of \cite{onedim}. There the authors prove that for a subanalytic subset $X$ of $\mathbb{Z}_p^{m+1}$ there is a semialgebraic set $X'\subset \mathbb{Z}_p^{m'+1}$ such that for each $\pb\in\mathbb{Z}_p$ there is a $\qb\in\mathbb{Z}_p^{m'+1}$ with $X(\pb)=X'(\qb)$. In their context{\em subanalytic} is in the sense of \cite{denefvdd} and semialgebraic is in Macintyre's Language. By the Quantifier Elimination theorem of \cite{denefvdd} one only needs to consider quantifier free definable subset of $\mathbb{Z}_p^{m+1}$ and one follows the argument of the proofs of \ref{thmA} above and Basic Lemma (1.2) of \cite{denefvdd} to get $y$ --the variable which corresponds to the fiber space-- to appear polynomially in each term of the formula that defines $X$ at the expense of introducing new variables that corresponds to terms that do not involve $y$.
\end{rem}


\end{document}